Theo van den Bogaart, august 2008

# The de Rham comparison theorem for Deligne-Mumford stacks

Let $K$ be a complete discrete valuation field of characteristic 0 with perfect residue field of characteristic $p \neq 0$, and let $\overline{K}$ be an algebraic closure of $K$. The de Rham comparison theorem in $p$-adic Hodge theory compares the $p$-adic étale and the de Rham cohomology of certain classes of $K$-schemes. Proofs are due to Faltings ([Fa]) and Tsuji ([Ts02]), although various persons have extended their results to larger classes of schemes. The present chapter is concerned with a version for smooth and proper Deligne-Mumford stacks. This is the following statement—the terminology that is used will be explained in §§1–2.

**0.1. Statement.** *Let $\mathscr{X}$ be a proper smooth Deligne-Mumford stack over $K$, and let $m \in \mathbb{Z}$. Then $\mathrm{H}^m(\mathscr{X}_{\overline{K},\text{ét}}, \mathbb{Q}_p)$ is a de Rham representation and there is a natural isomorphism*

$$D_{\mathrm{dR}}(\mathrm{H}^m(\mathscr{X}_{\overline{K},\text{ét}}, \mathbb{Q}_p)) \simeq \mathrm{H}^m_{\mathrm{dR}}(\mathscr{X}/K).$$

The main theorem of [B-E] fully describes the $p$-adic étale cohomology of a particular class of smooth and proper Deligne-Mumford stacks defined over $\mathbb{Q}_p$. In the last section of this article, a similar result is obtained for the Hodge structure of the associated analytic stack. This result would be an immediate corollary of the main theorem, if the de Rham comparison theorem in $p$-adic Hodge theory would be valid for smooth and proper Deligne-Mumford stacks. This is the motivation for the present chapter.

There are five parts. The first one (§§1–2) recalls certain facts about categories and functors that appear in Fontaine theory; after that it deals with the definition of $\ell$-adic étale and de Rham cohomology for Deligne-Mumford stacks.

The other parts are independent of each other and can be read in any order. Part 2 (§§3–4) describes a method, based on an idea by T. Tsuji and T. Saito, to extend the comparison theorem from schemes to stacks under the assumption that both étale and de Rham cohomology satisfy the formal properties of a Weil cohomology. These properties are studied in the third part (§§5–7), inspired by an approach due to Lafforgue. Unfortunately, the issue cannot be fully resolved.

Part 4 (§8) gives a proof of the comparison theorem, based on certain comparison theorems for simplicial schemes by Kisin or Tsuji. The last part (§§9–10) are the appendices, which contains certain technical results of independent interest.

We will now give an overview of the ideas in this chapter.

*The comparison theorem*

Both approaches towards an extension of the comparison theorem to Deligne-Mumford stacks rely on the existence of certain schemes whose cohomology can be related to that of the stack. The main existence result is a version of Chow's lemma for Deligne-Mumford stacks, first treated by Vistoli and generalized by Laumon and Moret-Bailly to the following result.

**0.2. Theorem** ([L-MB, 16.6.1]). *Let $\mathscr{X}$ be a Deligne-Mumford stack of finite type over a noetherian scheme $S$. There exists a quasi-projective $S$-scheme $X$ together with a proper, surjective and generically étale $S$-morphism $X \to \mathscr{X}$.* □

The first strategy for a proof of the comparison theorem, learned from Tsuji (and Kato, see [Ts02, proof of Prop. A5]), is to cut out the cohomology of a stack inside the cohomology of a scheme. The argument is based entirely on the formalism of correspondences. Let us sketch the argument, which is very easy to explain. Precise arguments and references appear in the sections below.

Let $\mathscr{X}$ be a connected, proper and smooth Deligne-Mumford stack over $K$, and let $m \in \mathbb{Z}$. The étale cohomology space $\mathrm{H}^m_{\text{ét}}(\mathscr{X}) := \mathrm{H}^m(\mathscr{X}_{\overline{K},\text{ét}}, \mathbb{Q}_p)$ is a $p$-adic representation of the absolute Galois group of $K$. The Fontaine functor $D_{\mathrm{dR}}$ associates to such a representation a $K$-vector space equipped with a decreasing filtration. The de Rham theorem for $\mathscr{X}$ that we want to prove says that $D_{\mathrm{dR}}(\mathrm{H}^i_{\text{ét}}(\mathscr{X}))$ is naturally isomorphic to the de Rham cohomology $\mathrm{H}^i_{\mathrm{dR}}(\mathscr{X}/K)$ of $\mathscr{X}$.

Now by the above version of Chow's lemma, combined with resolution of singularities, there exists a generically finite cover $f: Y \to \mathscr{X}$ of $\mathscr{X}$ by a smooth and proper $K$-scheme $Y$. Formal cohomological machinery then shows that $f_* \circ f^*: \mathrm{H}^i(\mathscr{X}) \to \mathrm{H}^i(\mathscr{X})$ is multiplication by a positive integer, where $\mathrm{H}^i(\mathscr{X})$ can be either étale or de Rham cohomology. So $f^*$ is injective and $f_*$ is surjective and in particular $f^*$ realizes $\mathrm{H}^i(\mathscr{X})$ as a subspace of $\mathrm{H}^i(Y)$. Faltings' theorem provides a functorial isomorphism $D_{\mathrm{dR}}(\mathrm{H}^i_{\text{ét}}(Y)) \xrightarrow{\sim} \mathrm{H}^i_{\mathrm{dR}}(Y)$ and it remains to see that this isomorphism preserves the image of $f^*$. But the image of $f^*$ equals the image of $f^* \circ f_*$ as $f_*$ is surjective. Again by formal cohomological machinery, the composition $f^* \circ f_*$ equals the correspondence given by a particular cycle in $Y \times_K Y$; thus we have removed all reference to the stack $\mathscr{X}$ from the problem. The comparison isomorphism for $\mathscr{X}$ now follows from the fact that the comparison isomorphism for schemes is compatible with correspondences.

One observes that the argument is purely formal. We have tried to develop a formalism that is the proper context for the above argument.

*Weil cohomology*

The argument above refers to 'formal cohomological machinery', a phrase that will be made exact below. In essence, we need that $\ell$-adic étale cohomology of smooth and proper Deligne-Mumford stacks over a field behaves 'just as for schemes': it satisfies Poincaré duality, the Künneth formula and there is a cycle map having the right compatibility properties. One recognizes here all the axioms for a Weil cohomology ([Kl68]).

That étale cohomology for Deligne-Mumford stacks satisfies these axioms may seem totally unsurprising. In fact, in the literature it is often taken for granted (for example, in the example in [B-M, §8]). However, there seems to be no proof of this fact in the literature. We provide proofs of some facts below, but there will remain one open problem: a proof that the intersection product of cycles is related to cup product of the corresponding classes in cohomology under the cycle map

$$A^r(X) \xrightarrow{\mathrm{Cl}} \mathrm{H}^{2r}_{\text{ét}}(X, \mathbb{Q}_\ell)(r) \qquad (r \in \mathbb{Z}).$$

Even in the case of proper, but not necessarily projective varieties, there seems to be no proof of this fact in the literature (see the remarks in section 7 below). The two main references, SGA $4\frac{1}{2}$ and [D-V], both rely on Chow's moving lemma, thus requiring the scheme to be embeddable in projective space.

Anyway, let us consider the positive results that we obtain. The theory of algebraic stacks, and especially its cohomology, was until recently a bit obscure. For a long time, the only real reference in the spirit of SGA was (the preliminary version of) the book by Laumon and Moret-Bailly [L-MB]. This book treats in great generality and with great precision the beginnings of the theory, but it does not touch many of the above mentioned subjects (see in particular [ibid., 18.7]). Around 2006, great progress is made by Laszlo and Olsson [L-O1], [L-O2], [L-O3], who develop a formalism of Grothendieck's six operations for a very general class of stacks (including Poincaré duality and the Künneth isomorphism). (They also correct a fundamental mistake made in [L-MB] concerning the cohomology of Artin stacks, indicated in [Beh03, Warning 5.3.12].)

We will make no use of the results of Laszlo and Olsson about stack cohomology here. Our treatment involves a comparison between the cohomology of stacks and their associated coarse moduli spaces. This has of course serious drawbacks: we have severe restrictions on the class of stacks and it works only for cohomology having a $\mathbb{Q}$-algebra as coefficients. The advantage is that it is down-to-earth, linking the cohomology of stacks to the cohomology of the more familiar schemes (or at least algebraic spaces).

Our treatment of cohomology is a variant of Lafforgue's, exhibited in the appendix to the article [Laf] in which he proved the Langlands correspondence for function fields (written before the work of Behrend or Laszlo and Olsson). Using comparison with coarse moduli spaces, Lafforgue developed part of the formalism of cohomology for a particular class of stacks, baptized serene stacks; which was sufficient for the applications he had in mind. The important property of such a serene stack is that Zariski locally it has a 'nice' (see [ibid., Déf. A.1]) cover. We will show how the constructions of Lafforgue can be adapted to Deligne-Mumford stacks. The serene property will be replaced by the property that a Deligne-Mumford stack is a quotient stack étale locally on the coarse moduli space.

Although we deal only with étale cohomology in §§5–7, we remark here that the proof of Weil's axioms for de Rham cohomology of a smooth and proper Deligne-Mumford stack defined over a field of characteristic 0 is easier—as is the case for schemes (cf., [Ha75] and [Ha70]).

*Simplicial methods*

Since the question whether étale cohomology is a Weil cohomology is not fully resolved, the method described above cannot be used to conclude that the comparison theorem holds (although we can conclude for instance that étale cohomology of a smooth and proper Deligne-Mumford $K$-stack is a so called de Rham representation, or that the comparison theorem holds for stacks admitting a 'nice' proper covering by a smooth scheme). Therefore, in §8 we take an entirely different approach.

The cohomology of a stack can be 'approximated' using simplicial schemes. This process must be regarded as some huge generalisation of Čech cohomology. Recently, comparison theorems for certain classes of simplicial schemes were established by Kisin ([Ki]) and (using different methods) by Tsuji ([Ts02]). In fact, Kisin shows how to derive the comparison theorem for (not necessarily smooth or proper) schemes from his results about simplicial schemes. The same procedure, with only very minor adjustments, can be applied to Deligne-Mumford stacks. The details of this construction are described in §8.

Let us end by remarking that another proof of the comparison theorem would be to use the coarse moduli space $\mathbb{X}$ associated to $\mathscr{X}$. If we work with $\mathbb{Q}$-algebras as coefficients, then the $\ell$-adic and de Rham cohomology spaces of the stack and it coarse moduli space coincide. The de Rham comparison theorem would therefore follow from comparison results for proper, but possibly singular schemes (see for example [Ki] or [Ya]).

*Notations*

The pair $(k, \ell)$ will always be used for a prime $\ell$ and a field $k$ of characteristic different from $\ell$. We fix a separable closure $\overline{k}$ of $k$ and denote by $G_k$ the Galois group of $\overline{k}/k$.

The ring of integers of any complete discrete valuation field $K'$ is denoted by $\mathcal{O}_{K'}$.

We fix a complete discrete valuation field $K$ of characteristic 0 with perfect residue field of characteristic $p \neq 0$, and an algebraic closure $\overline{K}$ of $K$. We will sometimes specialize to the case $k = K$ and $\ell = p$.

## 1. Fontaine theory

In this section, we recall some theory concerning Fontaine's de Rham functor. The main reference is [Fon82].

*Tensor categories*

Roughly speaking, in the definition we will adopt a tensor category is a category $\mathcal{T}$ in which the hom sets are vector spaces and for which a tensor product of objects is given, such that these structures satisfy "all nice properties one can imagine" (to quote Fontaine). The prototypical example is the category of finite dimensional vector spaces. To be precise, we adopt the terminology of [D-Mi]. Let $F$ be a field. A *tensor category over $F$* is an $F$-linear rigid tensor category such that $F \simeq \mathrm{End}(\underline{1})$ ([ibid., Def. 1.1, 1.7 and 1.15]). A *tensor functor* is an $F$-linear functor that preserves these structures ([ibid., Def. 1.8]). Note that the isomorphism $F \simeq \mathrm{End}(\underline{1})$ is unique, being an isomorphism of $F$-algebras. (Another often used reference is [Sa]. What we call a tensor category over $F$ is there called an $F$-linear rigid ACU $\otimes$-category such that $\mathrm{End}(\underline{1}) \simeq F$.) We will now treat two examples.

Let $\mathrm{Rep}_\ell(G_k)$ be the category of linear, continuous representations of $G_k := \mathrm{Gal}(\overline{k}/k)$ onto finite dimensional $\mathbb{Q}_\ell$-vector spaces (equipped with the $\ell$-adic topology). It forms an abelian tensor category over $\mathbb{Q}_\ell$ (it is even Tannakian). All the structures of a tensor category derive from the familiar constructions on the vector spaces underlying the objects of $\mathrm{Rep}_\ell(G_k)$. A unit object is a one-dimensional vector space on which the Galois group acts trivially.

We need another category, denoted $\mathrm{Fil}_k$. By definition, the objects of $\mathrm{Fil}_k$ are finite dimensional $k$-vector spaces $V$ equipped with a descending filtration

$$V \supset \cdots \supset \mathrm{Fil}^r V \supset \mathrm{Fil}^{r+1} V \supset \cdots$$

which is exhaustive and separating (which means that $\mathrm{Fil}^r V = V$ for $r \ll 0$ and $\mathrm{Fil}^r V = 0$ for $r \gg 0$). A morphism in this category is a linear map $f\colon V \to W$ such that $f(\mathrm{Fil}^r V) \subset \mathrm{Fil}^r W$. Such a morphism is called *strict* if $f(\mathrm{Fil}^r V) = f(V) \cap \mathrm{Fil}^r W$. The category $\mathrm{Fil}_k$ has again an obvious $k$-linear structure and a tensor product, making it a tensor category over $k$ ([Fon82, 3.2], [Fon94, 3.4]). A unit object is a one-dimensional vector space $V$ with the filtration given by ($\mathrm{Fil}^r V = 0 \iff r > 0$). Since the category $\mathrm{Fil}_k$ is not abelian (see the next remark), the following definitions are not redundant. One defines a sequence of morphisms

$$0 \longrightarrow V' \longrightarrow V \longrightarrow V'' \longrightarrow 0$$

in $\mathrm{Fil}_k$ to be a *short exact sequence* if all maps are strict and if the sequence of vector spaces obtained by forgetting the filtration is exact. An additive functor between $\mathrm{Fil}_k$ and an abelian category is *exact* if it preserves short exact sequences.

**1.1. Remarks.** 1. To see that $\mathrm{Fil}_k$ is not abelian, let $f\colon V \to V'$ be a morphism in $\mathrm{Fil}_k$ that is not an isomorphism, but which is an isomorphism on the underlying vector spaces. (One could take the same spaces, but equip the first with a filtration that is strictly coarser than that of the second space.) Then the image and coimage of $f$ exist, but do not agree.

2. We do not obtain a category if we simply add the demand that all morphisms need to be strict. Strictness is not preserved by composition. For example, let $\Delta = (k \oplus k)/k(1,1)$. Consider the following two maps between spaces whose $\mathrm{Fil}^0$ and $\mathrm{Fil}^1$ are as indicated, and with $\mathrm{Fil}^r = \mathrm{Fil}^0$ for $r \leq 0$ and $\mathrm{Fil}^r = 0$ for $r > 1$:

$$\begin{array}{ccccccc} \mathrm{Fil}^0 & : & k \oplus 0 & & k \oplus k & & \Delta \\ \cup & & \cup & \longrightarrow & \cup & \longrightarrow & \cup\,; \\ \mathrm{Fil}^1 & : & 0 & & 0 \oplus k & & \Delta \end{array}$$

the first map is the inclusion and the second one the canonical projection. One checks that both maps are strict, while their composition is not.

*Fontaine's de Rham functor*

We specialize to the case $k = K$, $\ell = p$ (see the introduction). The Fontaine ring $B_{\mathrm{dR}}$ is constructed in [Fon82]. It is a $\overline{K}$-algebra (it is even a field) equipped with a filtration and an action of the Galois group $G_K = \mathrm{Gal}(\overline{K}/K)$. The de Rham functor

$$\mathrm{Rep}_p(G_K) \xrightarrow{\overline{D}_{\mathrm{dR}}} \mathrm{Fil}_K$$

is defined by mapping a representation $W$ to the filtered vector space $\overline{D}_{\mathrm{dR}}(W) = (W \otimes_{\mathbb{Q}_p} B_{\mathrm{dR}})^{G_K}$.

Reflecting the fact that we will extend the comparison theorem by using only a formal argument, we do not need to know the exact definition of the Fontaine functor, but we only require some properties that we will discuss now (see [Fon82, 3.10 and 3.11] for the proofs). An object $V$ of $\mathrm{Rep}_p(G_K)$ is called a *de Rham representation* if the dimensions of the vector spaces underlying $V$ and $\overline{D}_{\mathrm{dR}}(V)$ coincide. The de Rham representations form a full subcategory $\mathrm{Rep}_{\mathrm{dR}}(G_K)$ of $\mathrm{Rep}_p(G_K)$, closed under subobjects, quotients, tensor products and duals; equipped with these structures it forms again a tensor category over $\mathbb{Q}_p$ (it is even a Tannakian subcategory).

The restriction $D_{\mathrm{dR}}$ of $\overline{D}_{\mathrm{dR}}$ to $\mathrm{Rep}_{\mathrm{dR}}(G_K)$ is an exact and faithful tensor functor. Let $V$ be an object of $\mathrm{Fil}_K$ and let $W$ be a de Rham representation. Then $V$ is isomorphic to $D_{\mathrm{dR}}(W)$ if and only if there is an isomorphism of $\overline{K}$-vector spaces $W \otimes_{\mathbb{Q}_p} B_{\mathrm{dR}} \simeq V \otimes_K B_{\mathrm{dR}}$ that preserves the filtrations and the Galois actions.

## 2. Cohomology

In this section we introduce the two objects of main interest. Let $\mathscr{X}$ be a smooth, finite type Deligne-Mumford stack over $k$. Let $m$ and $i$ be integers. We define $\ell$-adic étale cohomology $\mathrm{H}^m(\mathscr{X}_{\bar{k},\text{ét}}, \mathbb{Q}_\ell(i))$, which is an object of $\mathrm{Rep}_\ell(G_k)$, and de Rham cohomology $\mathrm{H}^m_{\mathrm{dR}}(\mathscr{X}/k, i)$, which is an object of $\mathrm{Fil}_k$. In both cases, the '$i$' refers to the $i$th Tate twist.

### *Functors from 2-categories*

Recall that the naive category whose objects are stacks and whose arrows are the usual morphisms of stacks is not very useful. For example, not all fibre products exist ([Gi, II.1.2.4.3]). The better thing to do is to consider stacks as living in a '2-category'. Roughly speaking, a 2-category is a category with an additional layer of morphisms (for the precise definition, see [Hak, I 1]). This means that for any two stacks $\mathscr{X}$ and $\mathscr{Y}$, the morphisms between them, called 1-morphisms, are themselves objects of a category $\mathrm{Hom}(\mathscr{X},\mathscr{Y})$. An arrow in this Hom-category is called a 2-morphism. Actually, the 2-category of stacks is a 2-groupoid: all 2-morphisms are isomorphisms. The 2-fibre product of stacks always exists (given by a square that commutes 'up to a 2-isomorphism').

Let $\mathscr{D}$ be a 2-category. A sub-2-category $\mathscr{V}$ of $\mathscr{D}$ will be called *full* if for any two objects $A$, $B$ of $\mathscr{V}$ the inclusion $\mathrm{Hom}_{\mathscr{V}}(A,B) \to \mathrm{Hom}_{\mathscr{D}}(A,B)$ between the categories of 1-morphisms $A \to B$ is an equivalence of categories.

If $\mathscr{D}$ is a 2-groupoid, define the category $\mathscr{D}^1$ having the same objects as $\mathscr{D}$, with as set of morphisms $\mathrm{Hom}_{\mathscr{D}^1}(A,B)$ the 2-isomorphism classes of 1-morphisms from $A$ to $B$ in $\mathscr{D}$, and with composition deduced from $\mathscr{D}$. Let $\mathscr{C}$ be an ordinary category. We define a *functor* $\mathscr{D} \to \mathscr{C}$ to be a functor $\mathscr{D}^1 \to \mathscr{C}$.

**2.1. Remarks.** 1. A 2-fibre product in $\mathscr{D}$ need not give an ordinary fibre product in $\mathscr{D}^1$. This fact is very inconvenient in the treatment of simpicial coverings in §8. Let us describe the most elementary counterexample.

Suppose $k$ is a separably closed field and put $S = \operatorname{Spec} k$. Let $\mathscr{D}$ be the 2-category of $k$-stacks. Let $G$ be a finite, non-trivial group. Consider $G$ as an $S$-group scheme that acts trivially on $S$. The fibre over $S$ of the quotient stack $[S/G]$ is a category that consists of only one object, corresponding to the $G$-torsor $G^t$ induced by $G$ itself, with automorphism group $G$: left multiplication by $g \in G$ defines a $G$-equivariant morphism $G^t \to G^t$. There is a unique morphism of $S$ to the quotient stack $[S/G]$. A 2-fibre product of $S \to [S/G]$ with itself is given by the diagram

$$\begin{array}{ccc} G & \longrightarrow & S \\ \downarrow & \overset{\epsilon}{\swarrow\!\!\!\!\swarrow} & \downarrow \\ S & \longrightarrow & [S/G] \end{array} .$$

Here $\epsilon$ is the 2-automorphism of $G \to [S/G]$ for which $\epsilon(g)$ is the automorphism of $G^t$ defined by $g$.

Now consider the commutative square

$$\begin{array}{ccc} S & \longrightarrow & S \\ \downarrow & & \downarrow \\ S & \longrightarrow & [S/G] \end{array}$$

in the 1-category $\mathscr{D}^1$ associated to $\mathscr{D}$. Each choice of an element $g \in G$ gives a pull-back of this diagram to a 2-commutative square

$$\begin{array}{ccc} S & \longrightarrow & S \\ \downarrow & \overset{\delta_g}{\swarrow\!\!\!\!\swarrow} & \downarrow \\ S & \longrightarrow & [S/G] \end{array}$$

in $\mathscr{D}$. The universal property of the 2-fibre product gives, for each choice of $\delta_g$, a morphism $S \to G$. Returning to $\mathscr{D}^1$, we conclude that although there certainly exists a map $S \to G$, it is not unique; hence $G$ is not the fibre product in $\mathscr{D}^1$.

2. Let $\mathscr{D}$ be a 2-groupoid and let $\mathscr{C}$ be an ordinary category. Let $\mathscr{C}^2$ be the 2-category having the same objects and 1-morphisms as $\mathscr{C}$, and whose only 2-morphisms are the identities. Denote by $\operatorname{Hom}(\mathscr{D}^1, \mathscr{C})$ the category whose objects are functors from $\mathscr{D}^1$ to $\mathscr{C}$, and whose arrows are the natural transformations. We can also consider *2-functors* from $\mathscr{D}$ to $\mathscr{C}^2$. They form the objects of a 2-category $\operatorname{Hom}(\mathscr{D}, \mathscr{C}^2)$; the arrows are the obvious 2-counterparts of natural transformations. (In general one has to differentiate between two notions of mappings between 2-categories, according to whether it is strictly associative or only associative up to a given 2-isomorphism. Also, the 2-functors actually form a 2-category in a natural way. We can neglect this: for 2-functors to $\mathscr{C}^2$ the notions of strong and weak 2-functors coincide, and the only 2-morphisms in $\operatorname{Hom}(\mathscr{D}, \mathscr{C}^2)$ would be the identities. For more information, see again [Hak].) There is a canonical functor

$$\operatorname{Hom}(\mathscr{D}, \mathscr{C}^2) \longrightarrow \operatorname{Hom}(\mathscr{D}^1, \mathscr{C})$$

which is an isomorphism of categories. Note that this statement has the form of an adjointness property. One could try to formalize this and generalize the above statement by, for example, include naturality considerations. The statement is not true if we replace $\mathscr{D}$ by a general 2-category. The right definition of $\mathscr{D}^1$ must then be that two 1-morphisms $f$ and $g$ are equivalent if they are linked by a chain of 2-morphisms.

*Cohomological formalism for stacks*

For the applications in this text we only need de Rham cohomology and $\ell$-adic étale cohomology of Deligne-Mumford stacks over a field. For some intermediate results it is convenient to have cohomology for general sheaves, and to work with derived functors. Therefore we will very briefly review the general construction of étale cohomology for stacks. The main references are [L-MB, §12, 13, 14 and 18] and [Laf, p. 205]. (As discussed in the introduction, the extensive work of Olsson and Laszlo on cohomology of arbitrary Artin stacks is not used here.)

Let $S$ be a quasi-separated scheme. See [D-Mu] or [L-MB] for the definition and main properties of Deligne-Mumford stacks. Following the conventions of [L-MB], all algebraic $S$-stacks are by definition quasi-separated. Let $\mathscr{X}$ be a Deligne-Mumford stack over $S$. Recall from [ibid., §12] the definition of the *étale site of* $\mathscr{X}$. This is a site whose objects are pairs $(U,u)$, called *étale opens*, with $U$ an algebraic $S$-space and with $u: U \to \mathscr{X}$ étale. A morphism $(U,u) \to (V,v)$ consists of a map $\varphi: U \to V$ together with a 2-isomorphism $u \xrightarrow{\sim} v\varphi$. The coverings of an object $(U,u)$ are the ones induced by the coverings of $U$ in the étale site of the algebraic space $U$. The category $\mathscr{X}_{\text{ét}}$ of sheaves on this sites is called the *étale topos of* $\mathscr{X}$.

If $f: \mathscr{X} \to \mathscr{Y}$ is a morphism of Deligne-Mumford $S$-stacks, then it induces a pair of adjoint functors $f^{-1}: \mathscr{Y}_{\text{ét}} \to \mathscr{X}_{\text{ét}}$ and $f_*: \mathscr{X}_{\text{ét}} \to \mathscr{Y}_{\text{ét}}$. We emphasize that the étale topology does not suffer from the mistake in [L-MB]: $f^{-1}$ is exact, hence the pair $(f^{-1}, f_*)$ defines a morphism of topoi. This can be proved in the classical way, the key fact being that for $f$ and $g$ composable morphisms of schemes, if $g$ and $gf$ are étale, then so is $f$ (which is false if we replace 'étale' by 'smooth'). To give some details: the problem is to show that $f^{-1}$ commutes with finite projective limits; by SGA 4 III 1.3 and I 5.4 this is true if the base change functor from étale $\mathscr{Y}$-schemes to étale $\mathscr{X}$-schemes defined by $f$ commutes with fibre products and difference kernels; for this one can immediately adapt the proof in [Mi, II.1.13].

We can now apply the formalism of derived categories. For instance, if $A$ is a ring of $\mathscr{X}_{\text{ét}}$, if $B$ is a ring of $\mathscr{Y}_{\text{ét}}$ and if $B \to f_*A$ is a given homomorphism, one has the derived functor $\mathrm{R}f_* : \mathrm{D}^+(\mathscr{X}_{\text{ét}}, A) \to \mathrm{D}^+(\mathscr{Y}_{\text{ét}}, B)$ from the derived category of bounded below complexes of $A$-modules to that of $B$-modules. We will encounter two specific situations later on: coherent sheaves and constructible $\ell$-adic sheaves.

If $i : \mathscr{Z} \to \mathscr{X}$ is a closed substack, and if $\mathscr{F}$ is a sheaf on the étale site of $\mathscr{X}$, one can form the subsheaf of sections with support in $\mathscr{Z}$, denoted $i_* i^! \mathscr{F}$ (see the appendix in § 9). This only depends on the closed subset induced by $\mathscr{Z}$ in the Zariski space $|\mathscr{X}|$ associated to $\mathscr{X}$. Thus we obtain cohomology with support in $\mathscr{Z}$: if $K \in \mathrm{ob}\,\mathrm{D}^+(\mathscr{X}_{\text{ét}}, A)$, put $\mathscr{H}_{\mathscr{Z}} K = i_* \mathrm{R} i^! K \in \mathrm{D}^+(\mathscr{X}_{\text{ét}}, A)$, and $\mathrm{R}_{\mathscr{Z}} f_* K = \mathrm{R} f_* \circ \mathscr{H}_{\mathscr{Z}} K \in \mathrm{ob}\,\mathrm{D}^+(\mathscr{Y}_{\text{ét}}, B)$. If $\mathscr{Z}'$ is a closed substack of $\mathscr{Z}$ there are functors $\mathscr{H}_{\mathscr{Z}'} \to \mathscr{H}_{\mathscr{Z}}$ and $\mathrm{R}_{\mathscr{Z}'} f_* \to \mathrm{R}_{\mathscr{Z}} f_*$.

Given $A$-modules $\mathscr{F}$, $\mathscr{G}$ and closed substacks $\mathscr{Z}$, $\mathscr{Z}'$ there are cup products

$$\mathrm{R}^m_{\mathscr{Z}} f_* \mathscr{F} \otimes_B \mathrm{R}^n_{\mathscr{Z}'} f_* \mathscr{G} \xrightarrow{\ \cup\ } \mathrm{R}^{m+n}_{\mathscr{Z} \cap \mathscr{Z}'} f_* (\mathscr{F} \otimes_A \mathscr{G}) \qquad (m, n \in \mathbb{Z}).$$

A construction (in the context of derived categories) is described in SGA $4\frac{1}{2}$-[*cycle*] 1.2. (But one can also adapt the construction in [Go, II 6.6].) This product is functorial, associative, anti-commutative and behaves well with respect to exact sequences. It commutes with the maps that extend the supports that are described in the last line of the previous paragraph.

*ℓ-Adic étale cohomology*

We specialize to the case that the base scheme $S$ is the spectrum of a regular ring of dimension $\leq 1$ on which $\ell$ is invertible. Suppose $f: \mathscr{X} \to \mathscr{Y}$ is a morphism between Deligne-Mumford $S$-stacks of finite type. In [L-MB, §18] the concept of a constructible sheaf is generalized to stacks, and it is shown that the derived functors of $f_*$ and $f^{-1}$ respect this notion. Working with projective systems of $\mathbb{Z}/\ell^n\mathbb{Z}$-modules (for varying $n$) modulo torsion (as in SGA 5 VI), the result is a pair of adjoint (when restricted to $\mathrm{D}^+$) functors

$$\mathrm{D}_c(\mathscr{Y}_{\text{ét}}, \mathbb{Q}_\ell) \xrightarrow{f^*} \mathrm{D}_c(\mathscr{X}_{\text{ét}}, \mathbb{Q}_\ell), \qquad \mathrm{D}^+_c(\mathscr{X}_{\text{ét}}, \mathbb{Q}_\ell) \xrightarrow{\mathrm{R}f_*} \mathrm{D}^+_c(\mathscr{Y}_{\text{ét}}, \mathbb{Q}_\ell)$$

between derived categories of constructible $\ell$-adic sheaves.

Let $\mathscr{X}$ be a finite type Deligne-Mumford $k$-stack. We are primarily interested in the $\mathbb{Q}_\ell$-vector spaces $\mathrm{H}^m(\mathscr{X}_{\bar{k},\text{ét}}, \mathbb{Q}_\ell(i))$ for various $i, m \in \mathbb{Z}$. Of course, these spaces vanish if $m < 0$. Let $s$ be the structure morphism of $\mathscr{X}$, and denote by $s_{\bar{k}}: \mathscr{X}_{\bar{k}} \to \operatorname{Spec} \bar{k}$ its base change to the separably closed field $\bar{k}$. Then $\mathrm{H}^m(\mathscr{X}_{\bar{k},\text{ét}}, \mathbb{Q}_\ell(i))$ equals $\mathrm{R}^m s_{\bar{k}*} \mathbb{Q}_\ell(i)$. This space is of finite dimension; a fact that follows from [L-MB, Cor. 18.3.2], combined with the standard technique to pass from finiteness results for cohomology with torsion coefficients to $\mathbb{Q}_\ell$-coefficients (as in the proof of [Mi, Lemma V.1.11]).

Given another finite type Deligne-Mumford $k$-stack $\mathscr{Y}$ with structure morphism $t$ and a $k$-morphism $f: \mathscr{X} \to \mathscr{Y}$, then by adjointness there is a natural map $\mathrm{R}t_{\bar{k}*} \mathbb{Q}_\ell(i) \to \mathrm{R}t_{\bar{k}*} \circ \mathrm{R}f_{\bar{k}*} \mathbb{Q}_\ell(i) = \mathrm{R}s_{\bar{k}*} \mathbb{Q}_\ell(i)$. Note that if $g$ is a morphism that is 2-isomorphic to $f$, then it induces the same maps on cohomology (the reason is that a 2-isomorphism $f \to g$ induces natural isomorphisms $f_* \to g_*$ and $f^{-1} \to g^{-1}$). So $\mathrm{H}^m(\mathscr{X}_{\bar{k},\text{ét}}, \mathbb{Q}_\ell(i))$ carries a Galois action and defines a contravariant functor, in the sense of the first subsection, from the 2-category of finite type Deligne-Mumford $k$-stacks to $\mathrm{Rep}_\ell(G_k)$. We will sometimes denote this functor by $\mathrm{H}^m_{\text{ét}}(-,i)$.

As explained above, multiplication on $\mathbb{Q}_\ell$ defines a cup product

$$\mathrm{H}^m(\mathcal{X}_{\bar{k},\text{ét}},\mathbb{Q}_\ell(i))\otimes_{\mathbb{Q}_\ell}\mathrm{H}^n(\mathcal{X}_{\bar{k},\text{ét}},\mathbb{Q}_\ell(j))\xrightarrow{\cup}\mathrm{H}^{m+n}(\mathcal{X}_{\bar{k},\text{ét}},\mathbb{Q}_\ell(i+j)).$$

*De Rham cohomology*

Let $\mathcal{X}$ be a smooth Deligne-Mumford $k$-stack. The structure sheaf $\mathcal{O}_\mathcal{X}$ of $\mathcal{X}$ is defined by $(\mathcal{O}_\mathcal{X})_{U,u}=\mathcal{O}_U$ and the obvious restriction maps ([L-MB, 12.7.1]). This gives rise to the de Rham complex of $\mathcal{O}_\mathcal{X}$ over $k$ (see [Il, II 1.1 and VIII 1.1, 2.1]). In fact, the module of differentials $\Omega^1_{\mathcal{X}/k}$ is locally free of finite rank and it is a coherent sheaf in the sense of [L-MB, §15], so the de Rham complex is the complex of sheaves of $k$-modules

$$\Omega^\bullet_{\mathcal{X}/k}\ :\quad 0\longrightarrow\mathcal{O}_\mathcal{X}\longrightarrow\Omega^1_{\mathcal{X}/k}\longrightarrow\Omega^2_{\mathcal{X}/k}\longrightarrow\cdots,$$

with $\Omega^r_{\mathcal{X}/k}=\bigwedge^r_{\mathcal{O}_\mathcal{X}}\Omega^1_{\mathcal{X}/k}$. We will now define a filtered complex $\Omega^\bullet_{\mathcal{X}/k,i}$ for each $i\in\mathbb{Z}$. All have $\Omega^\bullet_{\mathcal{X}/k}$ a their underlying complex. The filtration is defined by $\mathrm{Fil}^r\,\Omega^q_{\mathcal{X}/k,i}=0$ if $q<r+i$ and $\mathrm{Fil}^r\,\Omega^q_{\mathcal{X}/k,i}=\Omega^q_{\mathcal{X}/k}$ otherwise. Note that the wedge product defines a product $\Omega^\bullet_{\mathcal{X}/k,i}\otimes_{\mathcal{O}_\mathcal{X}}\Omega^\bullet_{\mathcal{X}/k,j}\to\Omega^\bullet_{\mathcal{X}/k,i+j}$ on filtered complexes.

De Rham cohomology is defined as the hypercohomology of these complexes:

$$\mathrm{H}^m_{\mathrm{dR}}(\mathcal{X}/k,i)=\mathrm{H}^m(\mathcal{X}_{\text{ét}},\Omega^\bullet_{\mathcal{X}/k,i})\qquad(m,i\in\mathbb{Z}).$$

The filtration on the de Rham complex induces a spectral sequence for hypercohomology, which defines a filtration on $\mathrm{H}^m_{\mathrm{dR}}(\mathcal{X}/k,i)$. If we forget the filtration, the spaces $\mathrm{H}^m_{\mathrm{dR}}(\mathcal{X}/k,i)$ for the different $i$ are all equal. From the spectral sequence, it also follows that if $\mathcal{X}$ is proper, then de Rham cohomology is a vector space of finite dimension (using finiteness of the cohomology of coherent sheaves, [L-MB, 15.6]). Thus de Rham cohomology defines a functor from the 2-category of proper and smooth Deligne-Mumford $k$-stacks to $\mathrm{Fil}_k$. The wedge product induces a product structure on these filtered cohomology spaces: one uses again the spectral sequence and the fact that the cup product respects exactness properties. (All this can be placed in the context of derived categories, but one has to introduce the concept of a 'derived category of filtered objects', because $\mathrm{Fil}_k$ is not abelian—see [Il]. In particular, this provides a treatment that avoids spectral sequences.)

*Abstract framework*

Inspired by [Jan, §6], we would like to place the above examples—étale and de Rham cohomology—in an abstract framework.

Let $(\mathcal{DM}/k)$ be the 2-category of finite type Deligne-Mumford $k$-stacks. Let $\mathcal{V}$ be a full sub-2-category of $(\mathcal{DM}/k)$ that is closed under isomorphisms, finite sums and finite products (all in the 2-categorical sense). Let $\mathcal{T}$ be a tensor category over a field $F$ with unit object $\underline{1}$.

**2.2. Definition.** A *cohomological structure* $\mathrm{H}$ *from* $\mathcal{V}$ *to* $\mathcal{T}$ consists of a collection of functors

$$(2.3) \qquad \mathrm{H}^m(-,i)\colon \mathcal{V}^{\mathrm{op}} \longrightarrow \mathcal{T} \qquad (m,i\in\mathbb{Z}),$$

together with natural transformations, called *cup products*,

$$(2.4) \qquad \mathrm{H}^m(-,i)\otimes\mathrm{H}^n(-,j) \xrightarrow{\cup} \mathrm{H}^{m+n}(-,i+j) \qquad (m,n,i,j\in\mathbb{Z})$$

that satisfy

(a) the maps $\cup$ are associative and $(-1,+1)$-graded commutative, which means that the morphism

$$\mathrm{H}^m(-,i)\otimes\mathrm{H}^n(-,j) \xrightarrow{\text{can.}} \mathrm{H}^n(-,j)\otimes\mathrm{H}^m(-,i) \xrightarrow{\cup} \mathrm{H}^{m+n}(-,i+j)$$

is $(-1)^{mn}$ times the map (2.4);

(b) there is a (necessarily unique) isomorphism $\alpha\colon\underline{1}\xrightarrow{\sim}\mathrm{H}^0(\operatorname{Spec}k,0)$ such that for any object $\mathcal{X}$ of $\mathcal{V}$ with structure morphism $s\colon\mathcal{X}\to\operatorname{Spec}k$, the following triangle commutes for all $m,i\in\mathbb{Z}$:

$$\begin{array}{ccc} \mathrm{H}^0(\mathcal{X},0)\otimes\mathrm{H}^m(\mathcal{X},i) & \xrightarrow{\ \cup\ } & \mathrm{H}^m(\mathcal{X},i)\,; \\ {\scriptstyle (\mathrm{H}^0(s,0)\circ\alpha)\otimes\mathrm{id}}\uparrow & \nearrow{\scriptstyle \text{can.}} & \\ \underline{1}\otimes\mathrm{H}^m(\mathcal{X},i) & & \end{array}$$

(c) given two objects $\mathcal{X}$ and $\mathcal{Y}$ in $\mathcal{V}$, the map

$$\mathrm{H}^m(\iota_{\mathcal{X}},i)+\mathrm{H}^m(\iota_{\mathcal{Y}},i)\ :\ \mathrm{H}^m(X\sqcup Y,i)\longrightarrow\mathrm{H}^m(X,i)\oplus\mathrm{H}^m(X,i),$$

induced by the inclusions $\iota_{\mathcal{X}}\colon\mathcal{X}\hookrightarrow\mathcal{X}\sqcup\mathcal{Y}$ and $\iota_{\mathcal{Y}}\colon\mathcal{Y}\hookrightarrow\mathcal{X}\sqcup\mathcal{Y}$, is an isomorphism;

(d) if $m$ is negative, $\mathrm{H}^m(\mathcal{X},i)=0$ for all $i$ and all objects $\mathcal{X}$ of $\mathcal{V}$.

If $f$ is a morphism in $\mathcal{V}$, its image under $\mathrm{H}^m(-,i)$ will often be denoted by $f^*$. We will use the same notation (and also its covariant version $f_*$) for the image of $f$ by some other functors we will encounter below, but we hope that this common abuse of notation will not cause any confusion.

**2.5. Remark.** If $\mathcal{T}$ is the category of finite dimensional $F$-vector spaces and one leaves out twisting, then one usually summarises (2.3), (2.4) and property (a) by considering a contravariant functor from $\mathcal{V}$ to the category of augmented, associative, graded, anti-commutative $F$-algebras (see e.g., [Kl68]). The analogous thing to do here would be to take 'the sum $\bigoplus_{m,i} \mathrm{H}^m(\mathcal{X},i)$' and say that this has the structure of an associative, graded 'algebra object' functorially in $\mathcal{X}$. Since countable sums need not exist in $\mathcal{T}$, this can only be done if we can extend $\mathcal{T}$ to a category where countable sums exist, and which is $F$-linear and has a tensor product, both compatible with the corresponding structures on $\mathcal{T}$. This is the case, for example, for the categories $\mathrm{Rep}_\ell(G_k)$ and $\mathrm{Fil}_F$.
We can also rephrase the definition as follows. Define $\mathcal{T}^{\mathbb{Z}\times\mathbb{Z}}$ as the category of functors from $\mathbb{Z}\times\mathbb{Z}$, regarded as a discrete category, to $\mathcal{T}$. So objects of $\mathcal{T}^{\mathbb{Z}\times\mathbb{Z}}$ are sequences $c=(c_{m,i})$ labeled by $m,i\in\mathbb{Z}$ and $\mathrm{Hom}(c,c')=\prod_{m,i}\mathrm{Hom}(c_{m,i},c'_{m,i})$. The category $\mathcal{T}^{\mathbb{Z}\times\mathbb{Z}}$ inherits an $F$-linear structure. For an object $c$ of $\mathcal{T}^{\mathbb{Z}\times\mathbb{Z}}$, consider a collection $\cup=(\cup_{m,i,n,j})$ of maps $\cup_{i,m,j,n}\colon c_{m,i}\otimes c_{n,j}\to c_{m+n,i+j}$ (for $m,i,n,j\in\mathbb{Z}$) such that associativity, and commutativity upto a factor $(-1)^{mn}$, hold. Such pairs $(c,\cup)$ form a category $\mathcal{T}^{\mathrm{alg}}$, where morphisms are those in $\mathcal{T}^{\mathbb{Z}\times\mathbb{Z}}$ that commute with the multiplication structures. A cohomological structure from $\mathcal{V}$ to $\mathcal{T}$ is then a functor $\mathcal{V}^{\mathrm{op}}\to\mathcal{T}^{\mathrm{alg}}$ that satisfies (b), (c) and (d) of the above definition.

## 3. Correspondences and Weil's axioms

We first discuss Vistoli's intersection theory for Deligne-Mumford stacks. Having this at our disposal, the construction of correspondences is only a formality. We then introduce axioms which allows one to carry correspondences over to cohomology.

### *Chow rings*

There are several different treatments in the literature concerning intersection theory for stacks. We will use the version of Vistoli, who develops the theory in the spirit of Fulton's book [Ful]. It should also be mentioned that Kresch [Kr] has greatly simplified part of the proofs of Vistoli. For us, the important result of Vistoli's theory is that, roughly speaking, one can define the rational Chow ring of a smooth, finite type Deligne-Mumford stack over a field that has the same look-and-feel as the Chow ring of a non-singular variety. We will now describe this in more detail. All proofs can be found in Vistoli's article [Vis].

Let $\mathcal{X}$ be a finite type Deligne-Mumford $k$-stack. For $d\geq 0$ and $(U,u)$ an étale open of $\mathcal{X}$, let $\mathrm{Z}^d(U,u)$ be the $\mathbb{Q}$-vector space generated by the integral closed subspaces of $U$ of codimension $d$ (without loss of generality, one may suppose $U$ to be a scheme). A morphism $(U,u)\to(V,v)$ defines a map $\mathrm{Z}^d(V,v)\to\mathrm{Z}^d(U,u)$ by pull-back of cycles ([Ful, 1.7]). In this way one obtains a presheaf $Z^d$ on the étale site of $\mathcal{X}$, which in fact is a sheaf. We put $Z^\bullet = \bigoplus_d Z^d$. The space of its global sections $Z^\bullet(\mathcal{X})$ is the graded vector space generated by the closed integral substacks of $\mathcal{X}$. Exploiting the sheaf properties, one can construct pull-backs along flat morphisms, push-forwards along proper representable morphisms and one can associate classes to arbitrary closed substacks of $\mathcal{X}$; it is also possible to construct a push-forward for non-representable proper morphisms.

There is also a sheaf $W^\bullet$: for each $d \geq 0$ and $(U, u)$ the vector space $W^d(U, u)$ is the direct sum of the multiplicative abelian group of non-zero rational functions on integral closed codimensions $d-1$ subspaces of $U$. The divisor map is denoted by $\delta: W^\bullet \to Z^\bullet$. The *(rational) Chow group* $A^\bullet(\mathcal{X})$ of $\mathcal{X}$ is defined as $Z^\bullet(\mathcal{X})/W^\bullet(\mathcal{X})$. Note that these are not the global sections of the sheaf cokernel of $\delta$, but only of its presheaf cokernel: rational equivalence is a global property. All previous constructions carry over: a morphism $f: \mathcal{X} \to \mathcal{Y}$ defines functorially a linear map $f^*: Z^d(\mathcal{Y}) \to Z^d(\mathcal{X})$ if $f$ is flat, and a map $f_*: Z^\bullet(\mathcal{X}) \to Z^\bullet(\mathcal{Y})$ if $f$ is proper (in fact, $f_*$ raises the degree by $\dim \mathcal{Y} - \dim \mathcal{X}$ if $\mathcal{X}$ and $\mathcal{Y}$ are equidimensional). A closed substack $\mathcal{Z}$ of codimension $d$ defines a class $[\mathcal{Z}] \in A^d(\mathcal{X})$. The construction of push-forwards for proper maps commutes with flat base change. Taking products of closed subschemes, one can form the exterior product $\times: A^d(\mathcal{X}) \otimes A^e(\mathcal{Y}) \to A^{d+e}(\mathcal{X} \times_k \mathcal{Y})$.

A representable morphism $f: \mathcal{X} \to \mathcal{Y}$ between finite type Deligne-Mumford $k$-stacks is a *regular local embedding of codimension d* if there exists a 2-commutative diagram

$$\begin{array}{ccc} X' & \xrightarrow{r} & Y \\ \downarrow & & \downarrow \\ \mathcal{X} & \xrightarrow{f} & \mathcal{Y} \end{array}$$

with both vertical maps surjective étale and with $r$ a regular closed embedding of codimension $d$ between schemes. (Let us emphasize that this square need not be cartesian: just take for $f$ a trivial covering of degree 2 between schemes, let the right vertical map equal $f$ and let the other two maps be the identities. In fact, the diagonal morphism of a smooth, pure-dimensional stack is a regular local embedding, but we cannot find a diagram as above that is cartesian unless the stack is representable by an algebraic space.) Like regular embeddings of schemes, this notion is stable under flat base change, although not for arbitrary base change. For such an $f$, Vistoli defines *Gysin maps*: for any finite type morphism of stacks $\mathcal{Y}' \to \mathcal{Y}$ this is a map of graded groups $f^! = f^!_{\mathcal{Y}'}: A^\bullet(\mathcal{Y}') \to A^\bullet(\mathcal{X}')$, where $\mathcal{X}' = \mathcal{X} \times_{\mathcal{Y}} \mathcal{Y}'$. Gysin morphisms commute with flat pull-backs and proper push forwards for Cartesian squares. There is also the following compatibility (omitted in [Vis]):

**3.1. Lemma.** *In the above situation, suppose the map $f': \mathcal{X}' \to \mathcal{Y}'$ is also a regular local embedding of codimension $d$. Let $\mathcal{Y}'' \to \mathcal{Y}'$ be a morphism of finite type and $x \in A^\bullet(\mathcal{Y}'')$. Then $f'^!_{\mathcal{Y}''} x = f^!_{\mathcal{Y}''} x$.*

PROOF. (Compare with the proof of [Ful, Th. 6.2(c)]. This proof is the only place in this text where we use the actual definition of the Gysin maps.) It suffices to prove the statement when $\mathscr{Y}''$ is a purely $r$-dimensional scheme for some integer $r$ and $x = [\mathscr{Y}'']$. Put $X'' = \mathscr{X} \times_{\mathscr{Y}} \mathscr{Y}'' \simeq \mathscr{X}' \times_{\mathscr{Y}'} \mathscr{Y}''$. Let $N$ (resp. $N'$) be the normal bundle of $f$ (resp. $f'$) pulled back to $X''$. There is a closed immersion $N' \to N$ (see [Vis, bottom of page 624]), which is an isomorphisms as $N'$ and $N$ are vector bundles of the same rank. Hence $f'^! x = f^! x$, for both are obtained by intersecting the normal cone of $X'' \to \mathscr{Y}''$ with the zero-section inside $N = N'$. □

Let $f: \mathscr{X} \to \mathscr{Y}$ be a morphism of finite type Deligne-Mumford $k$-stacks. Suppose $\mathscr{Y}$ is smooth of dimension $y$. Then the *graph* $\Gamma_f: \mathscr{X} \to \mathscr{Y} \times_k \mathscr{X}$ (sic) is a regular local embedding of codimension $y$. Indeed, choose étale presentations $Y \to \mathscr{Y}$ and $X \to \mathscr{X} \times_{\mathscr{Y}} Y$, such that $X$ and $Y$ are separated, finite type $k$-schemes. There is a canonical map $X \to Y$. Its graph $X \to Y \times_k X$ is a regular closed embedding of codimension $y$, and it lies over $\Gamma_f$.

As a particular example, one can take $\mathscr{X} = \mathscr{Y}$ (smooth) and let $f$ be the identity; its graph is the diagonal $\Delta: \mathscr{X} \to \mathscr{X} \times_k \mathscr{X}$. The *intersection product* on $A^\bullet(\mathscr{X})$ is defined by $x \cdot y = \Delta^!(x \times y)$ for $x, y \in A^\bullet(\mathscr{X})$. This product gives $A^\bullet(\mathscr{X})$ the structure of a graded $\mathbb{Q}$-algebra that is associative and commutative. The product has all the properties, familiar for schemes, of [Ful, §8.1].

*Correspondences*

Let $\mathscr{V}$ be a full sub-2-category of the 2-category of smooth and proper Deligne-Mumford $k$-stacks, closed under finite sums and finite products.

We are going to define a $\mathbb{Q}$-linear category $\mathscr{C}_{\mathscr{V}}$. It has the same objects as $\mathscr{V}$. The vector space of morphisms from $\mathscr{X}$ to $\mathscr{Y}$ is denoted $\mathrm{Corr}(\mathscr{X}, \mathscr{Y})$ and is called the space of *correspondences* from $\mathscr{X}$ to $\mathscr{Y}$. If one decomposes $\mathscr{X} = \sqcup \mathscr{X}_n$ into subspaces $\mathscr{X}_n$ of pure dimension $n$, then $\mathrm{Corr}(\mathscr{X}, \mathscr{Y}) = \sum A^n(\mathscr{X}_r \times \mathscr{Y})$ is a sum of rational Chow groups. Define the composition

$$\mathrm{Corr}(\mathscr{Y}, \mathscr{Z}) \otimes_{\mathbb{Q}} \mathrm{Corr}(\mathscr{X}, \mathscr{Y}) \longrightarrow \mathrm{Corr}(\mathscr{X}, \mathscr{Z})$$

by $(g, f) \mapsto g \circ f = (p_{13})_*(p_{12}^* f \cdot p_{23}^* g)$, where $p_{ij}$ denotes the projection from $\mathscr{X} \times \mathscr{Y} \times \mathscr{Z}$ onto the product of the $i$th and $j$th factor. That this composition is well-defined, associative, and that the diagonal serves as an identity, can be checked in the same way as the analogous results for varieties, because the formal properties of Chow rings for smooth stacks are the same as those for smooth varieties; references for these results concerning varieties are [Kl72, §4] or [Kl68, §1.3].

There is a functor

$$\mathscr{V}^{\mathrm{op}} \xrightarrow{h} \mathscr{C}_{\mathscr{V}}$$

which is the identity on objects and maps a morphism $f\colon\mathscr{X}\to\mathscr{Y}$ to the correspondence $f^* := (\Gamma_f)_*[\mathscr{X}] \in \mathrm{Corr}(\mathscr{Y},\mathscr{X})$ given by the graph of $f$. Now suppose that for each irreducible component $\mathscr{Y}'$ of $\mathscr{Y}$, each irreducible component of $f^{-1}\mathscr{Y}'$ has the same dimension as $\mathscr{Y}'$. Then the transpose $f_* := (\Gamma_f^t)_*[\mathscr{X}] \in \mathrm{Corr}(\mathscr{X},\mathscr{Y})$ of $f$ is defined by swapping the coordinates.

The *degree* $r$ of a separated dominant morphism $f$ of finite type between integral stacks is defined in [Vis, p. 620]. In general, $r$ is a non-negative rational number. If $f$ is representable, as it will be in all our applications of the following lemma, it is just the degree of the morphism of schemes obtained from a representation, so $r$ is an integer in that case. We extend the definition of degree to non-dominant morphisms by putting $r=0$ in that situation.

**3.2. Lemma.** *Let $\mathscr{M}$ and $\mathscr{N}$ be objects of $\mathscr{V}$. Suppose $\mathscr{M}$ and $\mathscr{N}$ are both connected and of the same dimension (recall that they are proper and smooth). Let $r\in\mathbb{Q}$ and suppose $f\colon\mathscr{M}\to\mathscr{N}$ is generically finite of degree $r$. Then $f_*\circ f^* = r\,\mathrm{id}^*_{\mathscr{N}}$.*

PROOF. Let $\Gamma_f\colon\mathscr{M}\to\mathscr{N}\times\mathscr{M}$ be the graph morphism of $f$, let $\Gamma_f^t$ be its transpose and let $\Delta_{\mathscr{N}}\colon\mathscr{N}\to\mathscr{N}\times\mathscr{N}$ be the diagonal. For $1\le i<j\le 3$, let $p_{ij}$ be the projection from $\mathscr{L}:=\mathscr{N}\times\mathscr{M}\times\mathscr{N}$ onto the product of the $i$th and $j$th factor. One has $p_{12}^*\Gamma_{f*}[\mathscr{M}] = (\Gamma_f\times\mathrm{id}_{\mathscr{N}})_*[\mathscr{M}\times\mathscr{N}]$ and $p_{23}^*\Gamma^t_{f*}[\mathscr{M}] = (\mathrm{id}_{\mathscr{N}}\times\Gamma_f^t)_*[\mathscr{N}\times\mathscr{M}]$. So writing out the definitions gives

$$f_*\circ f^* = p_{13*}\Big(((\Gamma_f\times\mathrm{id}_{\mathscr{N}})_*[\mathscr{M}\times\mathscr{N}])\cdot((\mathrm{id}_{\mathscr{N}}\times\Gamma_f^t)_*[\mathscr{N}\times\mathscr{M}])\Big).$$

To calculate the intersection product, consider the 2-cartesian square

$$\begin{array}{ccc} \mathscr{M} & \xrightarrow{\ \tilde{\delta}\ } & (\mathscr{M}\times\mathscr{N})\times(\mathscr{N}\times\mathscr{M}) \\ {\scriptstyle i}\big\downarrow & & \big\downarrow{\scriptstyle ((\Gamma_f\times\mathrm{id}_{\mathscr{N}})\times(\mathrm{id}_{\mathscr{N}}\times\Gamma_f^t))} \\ \mathscr{L} & \xrightarrow[\ \delta\ ]{} & \mathscr{L}\times\mathscr{L} \end{array}$$

with $\delta$ the diagonal morphism, $\tilde{\delta}=(\Gamma_f^t,\Gamma_f)$ and $i=(f,id,f)$. Choose étale presentations $N\to\mathscr{N}$ and $M\to\mathscr{M}\times_{\mathscr{N}}Y$ with $N$ and $M$ smooth, separated schemes of finite type over $k$. Let $f'\colon M\to N$ be the map lying over $f$. The map $M\to(M\times N)\times(N\times M)$ given by $x\mapsto(x,f(x),f(x),x)$ is a regular embedding of the same codimension as $\delta$; indeed, it is the graph of a morphism between smooth varieties. So $\tilde{\delta}$ is by definition a regular local immersion.

By [Vis, 5.5] $\tilde{\delta}^!([(\mathscr{M}\times\mathscr{N})\times(\mathscr{N}\times\mathscr{M})])=[\mathscr{M}]$. Hence the definition of the intersection product combined with Lemma 3.1 shows that $f_*\circ f^* = p_{13*}i_*[\mathscr{M}] = (f,f)_*[\mathscr{M}]$.

Now $(f,f)_*[\mathcal{M}] = \Delta_{\mathcal{N},*} f_*[\mathcal{M}]$. By definition, $f_*[\mathcal{M}] = r[\mathcal{N}]$ (recall that $r = 0$ if $f$ is not surjective). So $f_* \circ f^* = r\Delta_{\mathcal{N}}[\mathcal{N}] = r\,\mathrm{id}^*_{\mathcal{N}}$. □

*Weil cohomology*

Let $\mathcal{V}$ be a full sub-2-category of the 2-category of smooth and proper Deligne-Mumford $k$-stacks that is closed under isomorphisms, finite sums and finite products (all in the 2-categorical sense). Let $F$ be a field of characteristic 0 and let $\mathcal{T}$ be a tensor category over $F$.

**3.3. Definition.** A cohomological structure H from $\mathcal{V}$ to $\mathcal{T}$ (see Definition 2.2), is a *Weil cohomology* if the following axioms (K), (P), (C) and (P–C) are satisfied.

(K) (*Künneth isomorphism*) Let $\mathcal{X}$, $\mathcal{Y}$ be objects of $\mathcal{V}$ and denote by $p_{\mathcal{X}}$ and $p_{\mathcal{Y}}$ the canonical projections from $\mathcal{X} \times \mathcal{Y}$ onto the first and second factor, respectively. These projections induce maps $p^*_{\mathcal{X}}$ and $p^*_{\mathcal{Y}}$ on cohomology, whose cup products results in maps

$$\bigoplus_{h+l=m} \mathrm{H}^h(\mathcal{X}, i) \otimes \mathrm{H}^l(\mathcal{Y}, j) \longrightarrow \mathrm{H}^m(\mathcal{X} \times \mathcal{Y}, i+j) \qquad (i,j,m \in \mathbb{Z})$$

(note that the sum is finite by (d) in Definition 2.2). We require that these maps are isomorphisms.

(P) (*Poincaré duality*) Let $\mathcal{X}$ be an object of $\mathcal{V}$ and suppose that it is geometrically irreducible of dimension $x$. We require that $\mathrm{H}^{2x}(\mathcal{X}, x)$ is isomorphic to $\underline{1}$, the unit object of $\mathcal{T}$. We demand furthermore that for all $m, i \in \mathbb{Z}$ the morphism $\mathrm{H}^m(\mathcal{X}, i) \to \mathrm{H}^{2x-m}(\mathcal{X}, x-i)^\vee$ that is obtained from the composition

$$\mathrm{H}^m(\mathcal{X}, i) \otimes \mathrm{H}^{2x-m}(\mathcal{X}, x-i) \xrightarrow{\cup} \mathrm{H}^{2x}(\mathcal{X}, x) \xrightarrow{\sim} \underline{1}$$

is an isomorphism. (Note that this property does not depend on a particular choice for the isomorphism.)

(C) (*existence of cycle map*) For each $r \in \mathbb{Z}$ and each object $\mathcal{X}$ of $\mathcal{V}$ there exists a $\mathbb{Q}$-linear map

$$A^r(\mathcal{X}) \xrightarrow{\mathrm{Cl}^r_{\mathcal{X}}} \mathrm{Hom}(\underline{1}, \mathrm{H}^{2r}(\mathcal{X}, r)).$$

These *cycle maps* satisfy the following properties:

(i) they commute with pull-backs along flat morphisms;

(ii) they map intersection products to cup products;

(iii) $\mathrm{Cl}^0_{\mathrm{Spec}\,k}$ is not the zero map.

(P–C) (*compatibility*) It must be possible to choose for each irreducible $\mathcal{X}$ of dimension $x$ an isomorphism $\tau_{\mathcal{X}}: \mathrm{Hom}(\underline{1}, \mathrm{H}^{2x}(\mathcal{X}, x)) \xrightarrow{\sim} F$ (which exists by (P)) such that the following compatibility is satisfied. For each morphism $f:\mathscr{X} \to \mathscr{Y}$ between irreducible objects and for each $r \in \mathbb{Z}$, for each $z \in A^r(\mathscr{X})$ and for each $t \in \mathrm{Hom}(\underline{1}, \mathrm{H}^{2x-2r}(\mathscr{Y}, x-r))$

$$\tau_{\mathscr{Y}}(t \cup \mathrm{Cl}_{\mathscr{Y}} f_*(z)) = \tau_{\mathscr{X}}(f^* t \cup \mathrm{Cl}_{\mathscr{X}}(z)).$$

Let us make a few observations concerning this definition. First note that the cycle maps commute with sums and products. Indeed, if $\mathscr{U} \hookrightarrow \mathscr{X}$ and $\mathscr{V} \hookrightarrow \mathscr{Y}$ are closed integral substacks then

$$\mathrm{Cl}_{\mathscr{X}\times\mathscr{Y}}[\mathscr{U} \times \mathscr{V}] = \mathrm{Cl}_{\mathscr{X}\times\mathscr{Y}}(p_{\mathscr{X}}^*[\mathscr{U}] \cdot p_{\mathscr{Y}}^*[\mathscr{V}]) = p_{\mathscr{X}}^* \mathrm{Cl}_{\mathscr{X}}[\mathscr{U}] \cup p_{\mathscr{Y}}^* \mathrm{Cl}_{\mathscr{Y}}[\mathscr{V}]$$

and

$$(i_{\mathscr{X}}^* + i_{\mathscr{Y}}^*)\mathrm{Cl}_{\mathscr{X}\sqcup\mathscr{Y}}([\mathscr{U}] + [\mathscr{V}]) = \mathrm{Cl}_{\mathscr{X}}[\mathscr{U}] + \mathrm{Cl}_{\mathscr{Y}}[\mathscr{V}],$$

where $p_{\mathscr{X}}$, $p_{\mathscr{Y}}$ denote the canonical projections from $\mathscr{X} \times \mathscr{Y}$, and $i_{\mathscr{X}}$, $i_{\mathscr{Y}}$ the embeddings in $\mathscr{X} \sqcup \mathscr{Y}$ (which are all flat).

Since the cycle map commutes with cup products (Ciii) implies that the composition

$$\mathbb{Q} = A^0(\mathrm{Spec}\, k) \xrightarrow{\mathrm{Cl}} \mathrm{Hom}(\underline{1}, \mathrm{H}^0(\mathrm{Spec}\, k, 0)) = \mathrm{End}(\underline{1}) = F$$

is the unique inclusion of fields.

Poincaré duality makes it possible to define cohomological push-forwards. Indeed, for each irreducible $\mathscr{X}$ fix an isomorphism $\mathrm{H}^{2x}(\mathscr{X}, x) \tilde{\to} \underline{1}$ (in a moment, we will see there is a canonical choice). Let $f:\mathscr{X} \to \mathscr{Y}$ be a morphism in $\mathscr{V}$ and suppose $\mathscr{X}$ and $\mathscr{Y}$ are of pure dimensions $x$ and $y$, respectively. The map $f^* = \mathrm{H}^{2x-m}(f, x-i)$ has a transpose, living in $\mathrm{Hom}(\mathrm{H}^{2x-m}(\mathscr{X}, x-i)^\vee, \mathrm{H}^{2x-m}(\mathscr{Y}, x-i)^\vee)$. Applying Poincaré duality twice, we obtain a morphism $f_*:\mathrm{H}^m(\mathscr{X}, i) \to \mathrm{H}^{m+2(y-x)}(\mathscr{X}, i+y-x)$. This defines a push-forward in cohomology in a (covariant) functorial way. However, its definition depends on choices of isomorphisms $\mathrm{H}^{2x}(\mathscr{X}, x) \tilde{\to} \underline{1}$. We will now see that the class maps enable us to make a canonical choice.

In fact, the maps $\tau_{\mathscr{X}}$ in (P–C) determine such maps: take $\tau_{\mathscr{X}}^{-1}(1)$. We will always choose the isomorphisms in this way; axiom (P–C) then says that the cycle maps commute with the push-forwards. In fact, it suffices to assume that (P–C) holds for maps that are finite or flat, since $f_* = p_{\mathscr{Y}*} \circ \Gamma_{f*}$.

Given Cl, there is only one choice for $\tau_{\mathscr{X}}$ such that $\tau_{\mathrm{Spec}\, k}$ is the canonical isomorphism of (b) in Definition 2.2. In fact, take $z \in A^x(\mathscr{X})$ such that $s_*(z) \neq 0$ (for example $z = [P]$ for some integral zero-dimensional substack $P$). Then $\mathrm{Cl}(s_*(z)) = d\,\mathrm{Cl}[\mathrm{Spec}\, k]$ for some $d \in \mathbb{Q}^\times$, so $d \mapsto \mathrm{Cl}_{\mathscr{X}}^x(z)$ determines an isomorphism $F \tilde{\to} \mathrm{Hom}(\underline{1}, \mathrm{H}^{2x}(\mathscr{X}, x))$ of which $\tau_{\mathscr{X}}$ is the inverse.

**3.4. Remark.** Although we will not use it here, the class map also commutes with arbitrary pull-backs in the following sense. Varying the stack $\mathcal{X}$, one can view $A^r(\mathcal{X})$ as a contravariant functor from $\mathcal{V}$ to $\mathbb{Q}$-vector spaces. A map $f:\mathcal{X}\to\mathcal{Y}$ is mapped to $f^*$, which is defined by $f^*[V]=p_{\mathcal{X}*}\big(([V]\times[\mathcal{X}])\cdot\Gamma_{f*}[\mathcal{X}]\big)$. One easily checks that Cl now defines a natural transformation of functors $A^r(-)\to\mathrm{H}^{2r}(-,r)$.)

Let $c\in\mathrm{Hom}(\underline{1},\mathrm{H}^m(\mathcal{Y}\times\mathcal{X},i))$. Denote by $\cup c$ the image of $c\otimes\mathrm{id}$ under

$$\begin{array}{c}\mathrm{Hom}(\underline{1},\mathrm{H}^m(\mathcal{Y}\times\mathcal{X},i))\otimes\mathrm{Hom}(\mathrm{H}^n(\mathcal{Y}\times\mathcal{X},j),\mathrm{H}^n(\mathcal{Y}\times\mathcal{X},j))\\ \Big\downarrow{\scriptstyle\text{can.}}\\ \mathrm{Hom}(\mathrm{H}^n(\mathcal{Y}\times\mathcal{X},j),\mathrm{H}^m(\mathcal{Y}\times\mathcal{X},i)\otimes\mathrm{H}^n(\mathcal{Y}\times\mathcal{X},j))\\ \Big\downarrow{\scriptstyle\cup}\\ \mathrm{Hom}(\mathrm{H}^n(\mathcal{Y}\times\mathcal{X},j),\mathrm{H}^{m+n}(\mathcal{Y}\times\mathcal{X},i+j))\quad.\end{array}$$

Then we can define the morphism $\mathrm{corr}_{n,j}(c)$ as the composition

$$\mathrm{H}^n(\mathcal{Y},j)\xrightarrow{p_{\mathcal{Y}}^*}\mathrm{H}^n(\mathcal{Y}\times\mathcal{X},j)\xrightarrow{\cup c}\mathrm{H}^{n+m}(\mathcal{Y}\times\mathcal{X},j+i)\xrightarrow{p_{\mathcal{X},*}}\mathrm{H}^{n+m-2y}(\mathcal{X},j+i-y).$$

(We suppose a choice for cycle maps has been made.) The map corr is $F$-linear in $c$.

The following type of statement is well-known from e.g., [Kl68].

**3.5. Proposition.** *Suppose* H *is a Weil cohomology and fix the cycle maps. For all* $m,i\in\mathbb{Z}$, *there exists a unique functor* $\mathcal{C}_{\mathcal{V}}\to\mathcal{T}$ *such that the composition*

$$\mathcal{V}^{\mathrm{op}}\xrightarrow{h}\mathcal{C}_{\mathcal{V}}\longrightarrow\mathcal{T}$$

*equals* $\mathrm{H}^m(-,i)$ *and which maps* $c\in\mathrm{Hom}(h(\mathcal{Y}),h(\mathcal{X}))$, *for objects* $\mathcal{X}$ *and* $\mathcal{Y}$ *of* $\mathcal{V}$ *that are connected, to* $\mathrm{corr}_{m,i}(\mathrm{Cl}_{\mathcal{Y}\times\mathcal{X}}(c))\in\mathrm{Hom}(\mathrm{H}^m(\mathcal{Y},i),\mathrm{H}^m(\mathcal{X},i))$.

PROOF. Formally, the proof is identical to the proofs given in the literature, where $\mathcal{T}$ is the category of modules. One can consult the discussion concerning a map called '$T$' in [Kl72, §4]. □

Abusing notation, we denote the functor $\mathcal{C}_{\mathcal{V}}\to\mathcal{T}$ by $\mathrm{H}^m(-,i)$ again.

**3.6. Remark** (*continuation of Remark 2.5*). Recall that the functors $\mathrm{H}^m(-,i)$ for the different $m,i$ form a functor $\mathrm{H}:\mathcal{V}^{\mathrm{op}}\to\mathcal{T}^{\mathbb{Z}\times\mathbb{Z}}$. The proposition says that H induces a functor $\tilde{\mathrm{H}}:\mathcal{C}_{\mathcal{V}}\to\mathcal{T}^{\mathbb{Z}\times\mathbb{Z}}$.

One can easily extend the proposition if one adds morphisms of degree different from 0 as follows. First extend correspondences to graded objects by defining $\mathrm{Corr}^r(\mathcal{Y},\mathcal{X})=A^{y+r}(\mathcal{Y}\times\mathcal{X})$ for $\mathcal{Y}$ of pure dimension $y$. The composition of correspondences extends to a composition on these graded sets, thus forming a category $\mathcal{C}_{\mathcal{V}}^{\bullet}$ whose objects are those of $\mathcal{V}$ and with $\mathrm{Hom}(\mathscr{Y},\mathscr{X}) = \bigoplus_r \mathrm{Corr}^r(\mathscr{Y},\mathscr{X})$. Define $\mathscr{T}^\bullet$ as the category having the same objects as $\mathscr{C}^{\mathbb{Z}\times\mathbb{Z}}$, and with $\mathrm{Hom}(c,c') = \bigoplus_r \prod_{n,i} \mathrm{Hom}(c_{n,i}, c'_{n+2r,i+r})$. Then $\widetilde{\mathrm{H}}$ extends to a functor $\mathscr{C}^\bullet_{\mathscr{V}} \to \mathscr{T}^\bullet$ that respects the grading on the hom sets.
Suppose $\mathscr{T}$ is a pseudo-abelian category. Then one can exhibit motivic versions of the proposition and of the generalizations suggested in the above remark. In short, define the category of *effective (Chow) motives* $\mathscr{M}^+_{\mathscr{V}}$ as the pseudo-abelian envelope of $\mathscr{C}_{\mathscr{V}}$, i.e., the category obtained by formally adding the images of all idempotents (see [Kl72]). Then one may replace $\mathscr{C}_{\mathscr{V}}$ by $\mathscr{M}^+_{\mathscr{V}}$ in the proposition. One can extend in the usual way to the category of *full motives* $\mathscr{M}_{\mathscr{V}}$; see [D-Mi, §6], or [Sch]. One then obtains a functor

$$\mathscr{M}_{\mathscr{V}} \longrightarrow \mathscr{T}^\bullet.$$

A construction of Chow motives for Deligne-Mumford stacks can be found in [B-M].

*Projective schemes*

Let $(\mathscr{DM}^{\mathrm{sp}}/k)$ be the 2-category of smooth and proper Deligne-Mumford stacks over $k$. Let $(\mathrm{Sch}^{\mathrm{sP}}/k)$ be the subcategory of smooth and projective $k$-schemes. We have seen that étale and de Rham cohomology give cohomological structures $\mathrm{H}_{\text{ét}}$ and $\mathrm{H}_{\mathrm{dR}}$ from $(\mathscr{DM}^{\mathrm{sp}}/k)$ to $\mathrm{Rep}_\ell(G_k)$ or $\mathrm{Fil}_k$, respectively. Of course, we can restrict these structures to $(\mathrm{Sch}^{\mathrm{sP}}/k)$.

There is the following famous theorem.

**3.7. Theorem.** *Étale cohomology of smooth, projective $k$-schemes $\mathrm{H}_{\text{ét}}$ is a Weil cohomology.* □

Proofs of this theorem are scattered throughout SGA 4, SGA $4\frac{1}{2}$ and SGA 5. One can also consult [Mi]. There is a canonical choice for cycle maps, obtained from the construction of fundamental classes; see for example [Mi, VI 6].

**3.8. Theorem.** *If $k$ has characteristic zero, then de Rham cohomology of smooth, projective $k$-schemes $\mathrm{H}_{\mathrm{dR}}$ is a Weil cohomology.* □

A reference for this theorem is [Ha75]; a proof of Poincaré duality can also be found in [Ha70, III 8.5]. Since we have restricted attention to schemes over a field of characteristic zero, some parts of the theorem can also be derived using comparison results with ordinary singular cohomology. Again, fundamental classes give a canonical choice for the cycle maps.

The question whether étale cohomology of stacks leads to a Weil cohomology is the subject of the second part of this chapter, from section 5 onward.

## 4. The comparison theorem via correspondences

Recall that $K$ denotes a complete discrete valuation field of characteristic 0 with perfect residue field of characteristic $p \neq 0$. Tsuji's version of the de Rham comparison theorem immediately gives the following result.

**4.1. Theorem** (Tsuji [Ts02, Thm. A1]). *The Galois representation defined by the $p$-adic étale cohomology of a smooth and projective $K$-scheme is a de Rham representation. The following diagram is commutative up to a natural isomorphism*

$$\begin{array}{ccc} & \mathscr{C}_{(\mathrm{Sch}^{\mathrm{sP}}/K)} & \\ {}^{\mathrm{H}^m_{\text{ét}}(-,j)}\swarrow & & \searrow{}^{\mathrm{H}^m_{\mathrm{dR}}(-,j)} \\ \mathrm{Rep}_{\mathrm{dR}}(G_K) & \xrightarrow{\quad D_{\mathrm{dR}} \quad} & \mathrm{Fil}_K \end{array} .$$

*(We suppose that for both étale and de Rham cohomology the canonical choice for cycle maps has been made.)* □

**4.2. Remark** (*continuation of Remark 3.6*). One can make the generalizations suggested by Remark 3.6. The only thing to notice is that $\mathrm{Fil}_K$ is a pseudo-abelian category. Indeed, kernels exist in $\mathrm{Fil}_K$ (for all morphisms, not only projectors) and if $p\colon V \to V$ is a projector, then $\mathrm{Ker}(p) \oplus \mathrm{Ker}(1-p) \to V$ is an isomorphism, as this is true on the underlying vector spaces and the kernel is a subobject (hence strict). Therefore we obtain a diagram

$$\begin{array}{ccc} & \mathscr{M}(\mathrm{Sch}^{\mathrm{sP}}/K) & \\ {}^{\mathrm{H}_{\text{ét}}}\swarrow & & \searrow{}^{\mathrm{H}_{\mathrm{dR}}} \\ \mathrm{Rep}_{\mathrm{dR}}(G_K)^{\bullet} & \xrightarrow{\quad D^{\bullet}_{\mathrm{dR}} \quad} & \mathrm{Fil}^{\bullet}_K \end{array} ,$$

where $D^{\bullet}_{\mathrm{dR}}$ is the obvious extension of $D_{\mathrm{dR}}$ to the categories of graded objects.

**4.3. Theorem.** *Suppose that $\mathrm{H}_{\text{ét}}$ and $\mathrm{H}_{\mathrm{dR}}$ are Weil cohomologies on $(\mathscr{DM}^{\mathrm{sp}}/K)$. Then the de Rham comparison statement for smooth and proper Deligne-Mumford $K$-stacks (see Statement 0.1) is true.*

Proof. We may suppose that $\mathscr{X}$ is irreducible. Chow's lemma for Deligne-Mumford stacks (Theorem 0.2) says that there exists a proper, surjective, generically étale morphism from a projective scheme to $\mathscr{X}$. Combining this with resolution of singularities ([Hi]), there thus exists a surjective, generically finite morphism $f\colon Y \to \mathscr{X}$ where $Y$ is a projective and smooth scheme over $K$. We may assume that $Y$ is connected. Indeed, since $f$ proper, we can replace $Y$ by one of its connected components that maps surjectively onto the irreducible stack $\mathscr{X}$. Note that $Y$ and $\mathscr{X}$ have the same dimensions.

Choose cycle maps in such a way that they give the canonical choice for cycle maps when restricted to the category of schemes. The composition $f^* \circ f_*$ defines a morphism $g\colon h(Y) \to h(Y)$ in $\mathscr{C}_{(\mathrm{Sch}^{\mathrm{sp}}/K)}$. Now the comparison theorem for schemes (Theorem 4.1) gives a commutative diagram

$$\begin{array}{ccccc}
& & \xrightarrow{\;D_{\mathrm{dR}}(\mathrm{corr}_{m,0}(g))\;} & & \\
D_{\mathrm{dR}}\mathrm{H}^m_{\text{ét}}(Y,0) & \twoheadrightarrow & D_{\mathrm{dR}}\mathrm{H}^m_{\text{ét}}(\mathscr{X},0) & \hookrightarrow & D_{\mathrm{dR}}\mathrm{H}^m_{\text{ét}}(Y,0) \\
\downarrow \wr & & & & \downarrow \wr \\
\mathrm{H}^m_{\mathrm{dR}}(Y,0) & \twoheadrightarrow & \mathrm{H}^m_{\mathrm{dR}}(\mathscr{X},0) & \hookrightarrow & \mathrm{H}^m_{\mathrm{dR}}(Y,0) \\
& & \xrightarrow{\;\mathrm{corr}_{m,0}(g)\;} & &
\end{array}$$

where the two vertical maps are isomorphisms. Thus the restriction of the isomorphism at the right gives an isomorphism

$$D_{\mathrm{dR}}(\mathrm{H}^m_{\text{ét}}(\mathscr{X},0)) \xrightarrow{\sim} \mathrm{H}^m_{\mathrm{dR}}(\mathscr{X},0). \qquad \square$$

**4.4. Remark.** The statement that $\mathrm{H}^m_{\text{ét}}(\mathscr{X},\mathbb{Q}_p)$ is de Rham can be proved without the full set of assumptions made. In fact, it suffices to see that $f^*$ on étale cohomology is injective, as the property of being a de Rham representations is closed under subobjects. Injectivity follows from Poincaré duality and some elementary facts concerning trace maps and cohomology with compact support, which are all included in the next sections.

It is not hard to prove a more functorial version of the theorem, in the spirit of Theorem 4.1. For this we refer to the last page of Tsuji's article [Ts02, p. 368].

**4.5. Remarks.** Having the above comparison theorem at our disposal, we obtain a comparison theorem for algebraic spaces that are coarse moduli spaces for smooth and proper stacks. Such spaces are classified by Vistoli [Vis, 2.8]. One can speculate wether our method can be applied to other comparison theorems appearing in $p$-adic Hodge theory. For example, suppose the theory of crystalline cohomology can be extended from smooth schemes to smooth Deligne-Mumford stacks. (Although the author does not know any references for this fact, he does not see any reason why this cannot be done.) There exists a version of the crystalline comparison theorem similar to Theorem 4.1. The problem now is that we need more than resolution of singularities in characteristic 0, since we need a proper scheme which has a smooth model. Using De Jong's alteration, we would get a kind of 'potentially crystalline' comparison theorem. Another type of generalization would be towards the various kind of comparison results for the not-necessarily smooth or proper case. It is clear that a more elaborate cohomological formalism is necessary in that case, dealing (at least) also with cohomologies with compact support.

## 5. Coarse moduli spaces

The key to Lafforgue's treatment of cohomology of stacks lies in the fact that in favourable circumstances this cohomology can be identified with the cohomology of the coarse moduli space associated to the stack. The existence of a coarse moduli space was proven in [K-M] and we will begin this section by giving an overview of its main properties. After that, we state a lemma by Vistoli which is essential in the adaptation of Lafforgue's techniques to Deligne-Mumford stacks in the next sections. This lemma says that a Deligne-Mumford stack is a quotient stack for the action of a finite group on a scheme étale locally on the underlying coarse moduli space.

*Coarse moduli spaces: results of Keel and Mori*

Let $S$ be a quasi-separated, locally noetherian scheme. Let $\mathscr{X}$ be an algebraic $S$-stack. A *coarse moduli space of* $\mathscr{X}$ is an algebraic $S$-space $\mathbb{X}$ together with an $S$-morphism $q\colon \mathscr{X} \to \mathbb{X}$ that is universal for $S$-morphisms from $\mathscr{X}$ to an algebraic space. If such a space $\mathbb{X}$ exists, it is unique up to canonical isomorphism. One often says that $\mathbb{X}$ is *the* coarse moduli space of $\mathscr{X}$, with the map $q$ being understood.

We now recall a theorem by Keel and Mori, slightly adapted to our language.

**5.1. Theorem** (Keel and Mori [K-M]*). *Let $\mathscr{X}$ be a Deligne-Mumford $S$-stack that is separated and of finite type. Then a coarse moduli space $(\mathbb{X}, q)$ exists. Furthermore*

*(i) the algebraic space $\mathbb{X}$ is separated and of finite type over $S$;*

*(ii) for any flat morphism $Y \to \mathbb{X}$ of algebraic $S$-spaces, the space $Y$ together with the morphism $\mathscr{X} \times_{\mathbb{X}} Y \to Y$ obtained by base change is a coarse moduli space for $\mathscr{X} \times_{\mathbb{X}} Y$;*

*(iii) the map $q$ induces a bijection between $\mathbb{X}(F)$ and $\mathscr{X}_{\operatorname{Spec} F}/\simeq$ for every geometric point $\xi\colon \operatorname{Spec} F \to S$;*

*(iv) the morphism $q$ is surjective and universally open (this implies the often encountered statement that $q$ is universally submersive—see appendix 10).* □

**5.2. Corollary.**

*(v) If $\mathscr{X}$ is proper over $S$, then so is $\mathbb{X}$.*

*(vi) $\mathscr{X}$ is connected if and only if $\mathbb{X}$ is.*

---

* In the third sentence of [K-M], it is stated that all spaces must be assumed to be of finite type. This would imply that one has to add the finite type condition to the universality conditions in the definition of a coarse moduli space and in statement (ii). I think, also considering their proofs, that this is not what Keel and Mori intended. Anyhow, recently there have appeared generalizations of the theorem that do explicitly not require this assumption—nor even Noetherian assumptions for that matter. See [Co3] and [Rydh].

(vii) *$\mathscr{X}$ is irreducible if and only if $\mathbb{X}$ is.*

(viii) *$q\colon\mathscr{X} \to \mathbb{X}$ is locally quasi-finite.*

(ix) *Let $\mathscr{Y}$ be another separated, finite type Deligne-Mumford $S$-stack with coarse moduli space $r\colon\mathscr{Y} \to \mathbb{Y}$. Let $(Z,s)$ be the coarse moduli space of $\mathscr{X} \times_S \mathscr{Y}$. Then the canonical map $Z \to \mathbb{X} \times_S \mathbb{Y}$ is a universal homeomorphism.*

(x) *$\mathscr{X}$ and $\mathbb{X}$ have the same dimensions.*

PROOF. The map $q$ stays surjective after any base change $S' \to S$. Hence if $\mathscr{X} \to S$ is universally closed, so is $\mathbb{X} \to S$; this proves (v). Using the Zariski space associated to a stack, the 'only if' parts of (vi) and (vii) are elementary topological results. If $\mathscr{X}$ is a sum of two substacks, then by the universal property $\mathbb{X}$ is the sum of the coarse moduli spaces associated to these substacks; hence (vi). Suppose $\mathbb{X}$ is irreducible. If $\mathscr{U}_1$ and $\mathscr{U}_2$ are two non-empty open substacks of $\mathscr{X}$, then their images in the Zariski space $|\mathbb{X}|$ have a non-empty intersection $V$. Shrinking $\mathscr{U}_1$ and $\mathscr{U}_2$ if necessary, we may suppose $q(|\mathscr{U}_1|) = q(|\mathscr{U}_2|) \neq \emptyset$. But then (iii) implies $|\mathscr{U}_1| = |\mathscr{U}_2|$ and hence $\mathscr{U}_1 = \mathscr{U}_2$ by [L-MB, 5.4]. This proves (vii). For (ix), from the universal property one obtains the commutative diagram

$$\begin{array}{ccc} \mathscr{X} \times_S \mathscr{Y} & & \\ {\scriptstyle s}\downarrow & \searrow^{q\times r} & \\ Z & \xrightarrow[\alpha]{} & \mathbb{X} \times_S \mathbb{Y} \end{array} \quad .$$

Now $\alpha$ is surjective since $q \times r$ is, and $\alpha$ is universally injective since it is injective on the sets of geometric points (EGA I 3.5.5). Since $s$ is surjective and $q \times r$ is universally open, $\alpha$ is also universally open and hence a universal homeomorphism. Property (x) follows from (viii) and (iv): if $X \to \mathscr{X}$ is a presentation with $X$ quasi-compact, then $\dim(X) = \dim(\mathscr{X})$ and since the map $X \to \mathbb{X}$ is quasi-finite, surjective and open $\dim(X) = \dim(\mathbb{X})$ (EGA[IV 5.4.1]).

It remains to prove (viii). First note that 'locally quasi-finite' is a property of morphisms that is local on the source for the étale topology ([Kn, I 4.11]) and hence extends to arbitrary (not necessarily representable) morphisms of Deligne-Mumford stacks. Choose a presentation $f\colon X \to \mathscr{X}$ by a scheme $X$. Note that $qf$ is automatically schematic and locally of finite type. We may now reduce to the case that both $\mathbb{X}$ and $X$ are affine schemes. Let $x \in X$ and let $F$ be a separable closure of the residue field $k(qf(x))$. We have to show that $X_F := X \times_{\mathbb{X}} \operatorname{Spec} F$ has a finite number of points, or equivalently that $X_F(F)$ is finite. (Note that if $\overline{F}$ is an algebraic closure of $F$, then $X_{\overline{F}} \to X_F$ is a homeomorphism.) But if $x, y \in X_F(F)$ then $f(x) \simeq f(y)$ in $\mathscr{X}_F$ by (iii). Therefore, $\# X_F(F)$ is bounded by the number of points in an arbitrary fibre of the projection $X_F \times_{\mathscr{X}_F} X_F \to X_F$ to the first coordinate. But this projection is quasi-finite, being the base change of the étale morphism $f$. □

*A result of Vistoli*

Let us emphasize that neither the next lemma, nor the proof following it, are new. One can find the (almost similar) statement in [A-V, 2.2.3], or in [To]. Although [A-V] contains a sketch of a proof, a precise one is contained (as Toen notices) in the proof of a totally different statement in an article by Vistoli ([Vis, 2.8]). In view of the importance of the lemma for this text, and to have a clear cut proof, we give here Vistoli's proof with some details added.

**5.3. Lemma** (Vistoli). *Let $\mathscr{X}$ be a separated Deligne-Mumford $S$-stack of finite type and let $(\mathbb{X}, q)$ be its coarse moduli space. Let $x\colon \operatorname{Spec} K \to \mathbb{X}$ be an $S$-morphism with $K$ a field. There exists an étale morphism $g\colon U \to \mathbb{X}$ with $U$ an affine $S$-scheme, together with a lift $\operatorname{Spec} K \to U$ for which the following property holds: $U \times_{\mathbb{X}} \mathscr{X}$ is isomorphic to the quotient stack $[W/G]$ for some finite group $G$ acting on an affine $S$-scheme $W$.*

Proof. By [Kn, Thm. II 6.4] there exists an étale morphism $X' \to \mathbb{X}$ that factors $x$ with $X'$ a scheme. Therefore, we may assume that $\mathbb{X}$ is a scheme.

Consider triples $(U, u, y)$ that form a commutative triangle

$$\begin{array}{ccc} & & U \\ & \overset{y}{\nearrow} & \downarrow{\scriptstyle u} \\ \operatorname{Spec} K & \xrightarrow[x]{} & \mathbb{X} \end{array}$$

with $U$ an affine, connected scheme and $u$ étale. Given another such triple $(U', u', y')$ there exists at most one morphism $U \to U'$ that commutes with the other maps ([Mi, I.3.13]). So these triples define a partially ordered set $I$. The limit $U_\infty = \varprojlim U$ over the triples $(U, u, y) \in I$ exists, being a limit of affine schemes. It is isomorphic to the henselian scheme $\operatorname{Spec} \mathscr{O}^h_{\mathbb{X},x}$.

Choose an étale presentation $X \to \mathscr{X}$ with $X$ a separated $S$-scheme that is locally of finite type. For each $(U, u, y) \in I$, denote by $\mathscr{X}_U$ and $X_U$ the base change of respectively $\mathscr{X}$ and $X$ by $U$. Likewise, $\mathscr{X}_\infty$ and $X_\infty$ denote the base change by $U_\infty$ of $\mathscr{X}$ and $X$, respectively. This gives a diagram in the category of $S$-stacks with all squares 2-cartesian:

$$\begin{array}{ccccc} X_\infty & \longrightarrow & X_U & \longrightarrow & X \\ \downarrow & & \downarrow & & \downarrow \\ \mathscr{X}_\infty & \longrightarrow & \mathscr{X}_U & \longrightarrow & \mathscr{X} \\ \downarrow & & \downarrow & & \downarrow{\scriptstyle q} \\ U_\infty & \longrightarrow & U & \xrightarrow[u]{} & \mathbb{X} \end{array} \quad .$$

By EGA IV 8.2.5 we have $X_\infty = \varprojlim X_U$. We will now prove two successive claims.

Claim 1: *There is an affine connected component $W_\infty \subset X_\infty$ such that the restriction $W_\infty \to \mathcal{X}_\infty$ is surjective and finite.*

The map $U_\infty \to \mathbb{X}$ is flat (EGA IV 8.3.8ii) and hence $U_\infty$ is the coarse moduli space of $\mathcal{X}_\infty$ by (ii) of Theorem 5.1. Let $x' \in X_\infty$ be a closed point. As $U_\infty$ is henselian and $X_\infty \to U_\infty$ is separated, locally of finite type and quasi-finite at $x'$, by EGA IV 18.5.11(a⇒c) the scheme $W_\infty = \operatorname{Spec} \mathcal{O}_{X_\infty, x'}$ is a connected component of $X_\infty$ and the induced map $W_\infty \to U_\infty$ is finite. As $\mathcal{X}_\infty \to U_\infty$ is separated, this implies that $W_\infty \to \mathcal{X}_\infty$ is proper. But this map is also quasi-finite and hence $W_\infty \to \mathcal{X}_\infty$ is finite.

The fact that $W_\infty \to \mathcal{X}_\infty$ is an open (for it is étale) and closed morphism and that $\mathcal{X}_\infty$ is connected by (vi) of Corollary 5.2 implies that it is surjective. This proves claim 1.

Claim 2: *For some $(U, u, y) \in I$, there is an affine open subscheme $W \subset X_U$ such that $W \times_{\mathcal{X}_U} \mathcal{X}_\infty = W_\infty$ and such that the restriction $W \to \mathcal{X}_U$ is surjective, étale and finite.*

By EGA IV 8.2.11, there exists an open quasi-compact subscheme $W' \subset X_{U'}$ for some $(U', u, y) \in I$ such that $W_\infty = W' \times_{U'} U_\infty = W' \times_{\mathcal{X}_{U'}} \mathcal{X}_\infty$; the map $W_\infty \to U_\infty$ is affine, hence so is $W' \to U'$ (by [ibid., 8.10.5]); in particular $W'$ is affine. The map is clearly étale. It now suffices to show that $W' \times_{\mathcal{X}_{U'}} X_{U'} \to X_{U'}$ is surjective and finite after a base change $U \to U'$ obtained from some $(U', u', y') \in I$. This follows from EGA IV 8.10.5. Taking $W = W'_U$ gives claim 2.

From [L-MB, 6.1] the last claim implies that there exists a finite group $G$ acting on $W$ such that $\mathcal{X}_U \simeq [W/G]$. □

## 6. Comparison with the cohomology of the coarse moduli space

Let $G$ be a finite $S$-group that acts on an $S$-scheme $X$. Let $f: X \to \mathcal{Y}$ be a morphism to a Deligne-Mumford $S$-stack that is 2-invariant under the action (so the morphism $G \times_S X \to X$ is a morphism over $\mathcal{Y}$). If $\mathcal{F}$ is a sheaf on the étale site of $\mathcal{Y}$, then $f_* f^* \mathcal{F}$ carries a $G$-action and the canonical morphism $\mathcal{F} \to f_* f^* \mathcal{F}$ factors through the subsheaf of $G$-invariants $(f_* f^* \mathcal{F})^G$. If $G$ acts transitively on the fibres of $f$, then $\mathcal{F} \to (f_* f^* \mathcal{F})^G$ is an isomorphism.

Recall from §2 the construction of derived functors between derived categories of constructible $\ell$-adic sheaves on the étale sites of finite type Deligne-Mumford $k$-stacks. The next (well-known) lemma generalizes [B-E, Lemma 3.2]. It is our first example of an adaptation of a proof in [Laf] (in this case that of [ibid., A.3]) using Lemma 5.3 above.

**6.1. Lemma.** *Let $\mathcal{X}$ be a Deligne-Mumford $k$-stack that is separated and of finite type. Let $q: \mathcal{X} \to \mathbb{X}$ be the map to its coarse moduli space. Let $\mathcal{F}$ be a bounded-below complex of constructible $\ell$-adic sheaves on $\mathbb{X}$. Then the map $\mathcal{F} \to \mathrm{R}q_* q^* \mathcal{F}$, adjoint to the identity on $q^* \mathcal{F}$, is an isomorphism in $\mathrm{D}^+(\mathbb{X}_{\text{ét}}, \mathbb{Q}_\ell)$.*

PROOF. The question is étale local on $\mathbb{X}$. So by Lemma 5.3 we can assume that $\mathscr{X}$ is the quotient stack for a finite group $G$ acting on an affine, finite type $k$-scheme $W$. Denote by $f: W \to \mathscr{X}$ the quotient map; it is finite.

If we apply this to $q^*\mathscr{F}$ we obtain an isomorphism $\mathrm{R}q_*q^*\mathscr{F} \xrightarrow{\sim} \mathrm{R}q_*((\mathrm{R}f_*(qf)^*\mathscr{F})^G)$. Taking $G$-invariants commutes with $q_*$. Since $G$ is a finite group and $\#G$ is invertible in $\mathbb{Q}_\ell$, taking $G$-invariants is exact and hence derives trivially; in particular it commutes with $\mathrm{R}q_*$. From these two remarks it follows that $\mathrm{R}q_*((\mathrm{R}f_*(qf)^*\mathscr{F})^G) \simeq (\mathrm{R}(qf)_*(qf)^*\mathscr{F})^G$ which in turn is isomorphic to $\mathscr{F}$, the map $qf$ being finite (for it is proper and quasi-finite). □

This lemma tells us in particular that the cohomology of a stack and its coarse moduli space agree. In fact, if $f: \mathbb{X} \to S$ is any morphism of algebraic spaces then

$$\mathrm{R}^m(fq)_*q^*\mathscr{F} = \mathrm{R}^m f_*\mathrm{R}q_*q^*\mathscr{F} \simeq \mathrm{R}^m f_*\mathscr{F}. \tag{6.2}$$

There is in particular a canonical isomorphism $\mathrm{H}^m(\mathscr{X}_{\bar{k},\text{ét}}, \mathbb{Q}_\ell(i)) \simeq \mathrm{H}^m(\mathbb{X}_{\bar{k},\text{ét}}, \mathbb{Q}_\ell(i))$.

**6.3. Proposition.** *If $\mathscr{X}$ and $\mathscr{Y}$ are proper Deligne-Mumford $k$-stacks, then the canonical map (see Definition 3.3(K))*

$$\bigoplus_{h+l=m} \mathrm{H}^h(\mathscr{X}_{\text{ét}}, \mathbb{Q}_\ell(i)) \otimes_{\mathbb{Q}_\ell} \mathrm{H}^l(\mathscr{Y}_{\text{ét}}, \mathbb{Q}_\ell(j)) \longrightarrow \mathrm{H}^m((\mathscr{X} \times_k \mathscr{Y})_{\text{ét}}, \mathbb{Q}_\ell(i+j))$$

*is an isomorphism.*

PROOF. (See [Laf, A.8].) Let $\mathrm{Z}$ be the coarse moduli space associated to $\mathscr{X} \times_k \mathscr{Y}$. As we have seen (Corollary 5.2(ix)) the morphism $\alpha: \mathrm{Z} \to \mathbb{X} \times_k \mathbb{Y}$ is a universal homeomorphisms and as it is also of finite type, it is finite, surjective and radical (EGA IV 2.4.4). Therefore, the map $\mathrm{H}^n((\mathbb{X} \times_k \mathbb{Y})_{\text{ét}}, \mathbb{Q}_\ell(i+j)) \to \mathrm{H}^n(\mathrm{Z}_{\text{ét}}, \mathbb{Q}_\ell(i+j))$ induced by $\alpha$ is an isomorphism (SGA 4 VIII 1.2).

Since we have seen before the statement of the proposition that the cohomology of a stack agrees with that of its coarse moduli space, it suffices to show the proposition with $\mathscr{X}$ and $\mathscr{Y}$ replaced by $\mathbb{X}$ and $\mathbb{Y}$. This is the Künneth formula for algebraic spaces (for schemes, see e.g., [Mi, VI.8.5]). □

*Poincaré duality*

Recall (SGA 4 XVIII 3.1.4) that for a compactifiable morphism of schemes $f: X \to Y$ of finite type over $k$, the functor 'cohomology with compact support'

$$\mathrm{D}(X, \mathbb{Q}_\ell) \xrightarrow{\mathrm{R}f_!} \mathrm{D}(Y, \mathbb{Q}_\ell)$$

has a partial right adjoint

$$\mathrm{D}^+(X,\mathbb{Q}_\ell) \overset{\mathrm{R}f^!}{\longleftarrow} \mathrm{D}^+(Y,\mathbb{Q}_\ell).$$

(The notations $\mathrm{R}f_!$ and $\mathrm{R}f^!$ are misleading: in general they are not derived functors.) The definition of cohomology with compact support and this adjointsness result extend straightforwardly to algebraic spaces (this fact seems well-known and is frequently used in the literature, but the author does not know a straightforward reference; the results are of course a corollary of the much stronger and much more technical results of Laszlo and Olsson — [L-O1], [L-O2], [L-O3]).

**6.4. Theorem.** *Let $\mathcal{X}$ be a smooth open substack of pure dimension $x$ of a proper Deligne-Mumford $k$-stack (for some $x \in \mathbb{Z}$). Let $(\mathbb{X},q)$ be the coarse moduli space of $\mathcal{X}$ with structure morphism $s\colon \mathbb{X} \to \operatorname{Spec} k$. Then there is an isomorphism*

(6.5) $$\mathbb{Q}_\ell(x)[2x] \xrightarrow{\sim} \mathrm{R}s^!\mathbb{Q}_\ell$$

*in* $\mathrm{D}^+(\mathbb{X},\mathbb{Q}_\ell)$.

Proof. Note that $\mathbb{X}$ is compactifiable, so $\mathrm{R}s^!$ is well-defined. Let us first show that the presheaf on the étale site of $\mathbb{X}$ given by

(6.6) $$(U \xrightarrow{g} \mathbb{X}) \mapsto \operatorname{Hom}_{\mathrm{D}^+(U_{\text{ét}},\mathbb{Q}_\ell)}(\mathbb{Q}_\ell(x)[2x], g^*\mathrm{R}s^!\mathbb{Q}_\ell)$$

is a sheaf. By [BBD, Prop. 3.2.2] a sufficient condition for this is the vanishing of the sheaf $\underline{\operatorname{Ext}}^i(\mathbb{Q}_\ell(x)[2x],\mathrm{R}s^!\mathbb{Q}_\ell)$ for $i<0$. Since $\mathbb{X}$ has dimension $x$ we have $\mathrm{R}^j s^!\mathbb{Q}_\ell = 0$ for $j < -2x$ by SGA4[XVIII 3.1.7]. So there is a complex $I^\bullet$ of injective modules whose image in $\mathrm{D}(\mathbb{X}_{\text{ét}},\mathbb{Q}_\ell)$ is isomorphic to $\mathrm{R}s^!\mathbb{Q}_\ell$ such that $I^j=0$ for all $j<-2x$. By definition ([Ha66, II 3]) $\underline{\operatorname{Ext}}^i(\mathbb{Q}_\ell(x)[2x],\mathrm{R}s^!\mathbb{Q}_\ell)$ is the $i$-th cohomology of the complex

$$\prod_{j\in\mathbb{Z}} \underline{\operatorname{Hom}}(\mathbb{Q}_\ell(x)^{j+2x}, I^{j+\bullet})$$

(the $d$-maps are irrelevant). As the only non-trivial term in this product corresponds to $j=-2x$, the homology is zero for $i<0$; thus proving that (6.6) defines a sheaf.

Consider the étale site of $\mathbb{X}$ whose objects are affine schemes that are étale over $\mathbb{X}$. Let $C$ be the full subcategory consisting of objects $g\colon U \to \mathbb{X}$ for which $U \times_{\mathbb{X}} \mathcal{X}$ is isomorphic to the quotient stack for some finite group acting on an affine scheme, just as in Lemma 5.3. Then the lemma says that $C$ forms a sieve that covers $\mathbb{X}$. Each map $g$ is separated, étale and of finite type and therefore it is compactifiable. There is a localization isomorphism $g^*\mathrm{R}s^! \simeq \mathrm{R}(sg)^!$. The base change $\mathcal{X} \times_{\mathbb{X}} U$ is a quotient stack and [Laf, Prop. A.5] states that there is a natural isomorphism

$$\varphi_t : \mathbb{Q}_\ell(x)[2x] \xrightarrow{\sim} \mathrm{R}(sg)^!\mathbb{Q}_\ell \simeq g^*\mathrm{R}s^!\mathbb{Q}_\ell$$

in $\mathrm{D}^+(U_{\text{ét}},\mathbb{Q}_\ell)$. For different $U$ in $C$, these maps form a descent datum. Since (6.6) is a sheaf, the maps glue to a global isomorphism. □

**6.7. Remarks.** 1. An alternative proof of the theorem would be to copy verbatim the proof given in [Laf, Prop. A.5], excluding the localization step: $\mathcal{X} \times_{\mathbb{X}} \mathcal{U}$ is a smooth quotient stack of the right dimension.

2. If the ground field has characteristic zero, one can also prove directly that the coarse moduli space of a smooth stack satisfies Poincaré duality; this result is due to Satake. For details, see [B-E, Prop. 3.3].

If $\mathcal{X}$ is also proper, the adjoint of (6.5) gives a morphism

$$\mathrm{H}^{2x}(\mathcal{X}_{\bar{k},\text{ét}}, \mathbb{Q}_\ell(x)) \xrightarrow{\mathrm{Tr}_{\mathcal{X}}} \mathbb{Q}_\ell$$

called the *trace map*. One observes that it is an isomorphism if $\mathcal{X}$ is geometrically irreducible. In particular, the action of the Galois group $G_k$ on $\mathrm{H}^{2x}(\mathcal{X}_{\bar{k},\text{ét}}, \mathbb{Q}_\ell(x))$ is trivial.

As explained in SGA 4 XVIII 3.2.6 (where the 'RHom-version' of the adjointness property that is used there is obtained from [ibid., 3.1.10]), Poincaré duality (see Definition 3.3) follows from this theorem. In fact, let $\mathcal{X}$ be a smooth, proper Deligne-Mumford $k$-stack and let $t: \mathcal{X}_{\bar{k}} \to \operatorname{Spec} \bar{k}$ be the structure morphism of the base change of $\mathcal{X}$ to the separable closure of $k$. Then the adjointness of $\mathrm{R}t_*$ and $\mathrm{R}t^!$ combined with the above theorem gives a natural isomorphism

$$\mathrm{R}t_* \mathrm{R}\underline{\mathrm{Hom}}(\mathbb{Q}_\ell(i), \mathbb{Q}_\ell(x)[2x]) \xrightarrow{\sim} \mathrm{R}\mathrm{Hom}(\mathrm{R}t_* \mathbb{Q}_\ell(i), \mathbb{Q}_\ell),$$

whose $m$th homology is the desired isomorphism

$$\mathrm{H}^{2x-m}(\mathcal{X}_{\bar{k},\text{ét}}, \mathbb{Q}_\ell(x-i)) \xrightarrow{\sim} \mathrm{H}^m(\mathcal{X}_{\bar{k},\text{ét}}, \mathbb{Q}_\ell(i))^\vee.$$

It is advantageous to also have a notion of cohomology with compact support. This can be defined in the usual way, but instead we follow Lafforgue again and adopt an ad-hoc definition inspired by (6.2). Let $\mathcal{X}$ be an open substack of a proper Deligne-Mumford $k$-stack (we call $\mathcal{X}$ *compactifiable*). Define

$$\mathrm{H}_c^m(\mathcal{X}_{\bar{k}}, \mathbb{Q}_\ell(i)) = \mathrm{H}_c^m(\mathbb{X}_{\bar{k},\text{ét}}, \mathbb{Q}_\ell(i)),$$

where the object on the right is the cohomology with compact support of the coarse moduli space $\mathbb{X}_{\bar{k}}$. If $\mathcal{X}$ is also smooth of dimension $x$, the theorem above gives a trace map

$$\mathrm{H}_c^{2x}(\mathcal{X}_{\bar{k},\text{ét}}, \mathbb{Q}_\ell(x)) \xrightarrow{\mathrm{Tr}_{\mathcal{X}}} \mathbb{Q}_\ell.$$

Even better, this map exists if $\mathcal{X}$ is smooth on a dense open substack. For example, this is the case when $\mathcal{X}$ is reduced and $k$ is perfect.

## 7. The cycle map

Note that to give a morphism $\underline{1} \to \mathrm{H}^m(\mathscr{X}_{\bar{k},\text{ét}}, \mathbb{Q}_\ell(i))$ between objects of $\mathrm{Rep}_\ell(G_k)$ amounts to the same thing as giving an element of $\mathrm{H}^m(\mathscr{X}_{\bar{k},\text{ét}}, \mathbb{Q}_\ell(i))$ that is Galois-invariant. Hence for a connected, proper, smooth $k$-stack of dimension $x$, the inverse of the trace map

$$\mathrm{H}^{2x}(\mathscr{X}_{\bar{k},\text{ét}}, \mathbb{Q}_\ell(x)) \xrightarrow{\sim} \mathbb{Q}_\ell$$

provides an isomorphism $\tau_{\mathscr{X}}$ alluded to in (P–C) of Definition 3.3 and hence can be used to define push-forwards in cohomology. The goal of this section is to define cycle maps

$$A^r(\mathscr{X}) \xrightarrow{\mathrm{Cl}^r_{\mathscr{X}}} \mathrm{H}^{2r}(\mathscr{X}_{\bar{k},\text{ét}}, \mathbb{Q}_\ell(r))^{G_k} \qquad (r \in \mathbb{Z})$$

compatible with the maps $\tau_{\mathscr{X}}$.

Again, we follow Lafforgue ([Laf, §A2a] this time) in the definition of cycle maps. For $\mathscr{Z} \subset \mathscr{X}$ a closed substack of a Deligne-Mumford $k$-stack $\mathscr{X}$, recall from §2 that $\mathscr{H}_{\mathscr{Z}}$ is the derived functor 'sections with support in $\mathscr{Z}$'; in particular $\mathscr{H}^m_{\mathscr{Z}}\mathbb{Q}_\ell(i)$ is a sheaf on $\mathscr{X}$ for the étale topology, for each $m, i \in \mathbb{Z}$. The key to the construction of cycle classes is the following purity result, which is an immediate consequence of the analogous result for schemes (see SGA $4\frac{1}{2}$-[*cycle*] 2.2.8).

**7.1. Proposition** ([Laf, A.10]). *Let $c, i \in \mathbb{Z}$. Let $\mathscr{X}$ be a smooth finite type Deligne-Mumford $k$-stack and $\mathscr{Z} \hookrightarrow \mathscr{X}$ a closed substack of codimension $\geq c$. Then $\mathscr{H}^m_{\mathscr{Z}}\mathbb{Q}_\ell(i) = 0$ for all $m < 2c$.* □

(Notice that the cited proof of this theorem is valid for any ground field of characteristic different from $\ell$.)

The composite functor spectral sequence gives the following corollary to this proposition:

**7.2. Corollary.** $\mathrm{H}^{2c}_{\mathscr{Z}}(\mathscr{X}, \mathbb{Q}_\ell(c)) = \mathrm{H}^0(\mathscr{X}, \mathscr{H}^{2c}_{\mathscr{Z}}\mathbb{Q}_\ell(c))$. □

Thus, loosely speaking, '$U \mapsto \mathrm{H}^{2c}_{\mathscr{Z}}(U, \mathbb{Q}_\ell(c))$' is a sheaf on $\mathscr{X}$ so that we can define cycle classes locally and then glue. We will now describe this constructions in more detail.

First we briefly recall the constructions of cycle classes on schemes that is described in SGA $4\frac{1}{2}$-[*cycle*]. Let $X$ be a smooth scheme over the field $k$. The first step is to define the class of a divisor $D$: a divisor determines an element in the Picard group, which (using torsors) is isomorphic to $\mathrm{H}^1(X_{\text{ét}}, \mathbb{G}_m)$. The map from this space to $\mathrm{H}^2(X_{\text{ét}}, \mathbb{Q}_\ell(1))$ is provided by the Kummer sequence. Taking care of supports, one in fact obtains a class in the cohomology with support in $D$. Now suppose $Z \subset X$ is a local complete intersection of codimension $c$. Since locally $Z$ is the intersection of $c$ divisors, one can define $\mathcal{C}\ell(Z) \in \Gamma(X, \mathcal{H}_Z^{2c}\mathbb{Q}_\ell(c))$ locally by taking the cup product of the classes of these divisors, and then glue ([ibid., 2.1 and 2.2]). If $Z \subset X$ is an arbitrary integral closed subscheme of codimension $c$, it is a local complete intersection on some open piece whose complement has codimension $> c$. Another purity result then assures that the class of this open piece uniquely extends to a class $\mathcal{C}\ell(Z) \in \Gamma(X, \mathcal{H}_Z^{2c}\mathbb{Q}_\ell(c))$. By linearity, this construction extends to arbitrary cycles $z$ of codimension $c$ on $X$ and gives a class $\mathcal{C}\ell(z) \in \Gamma(X, \mathcal{H}_{|z|}^{2c}\mathbb{Q}_\ell(c))$, where $|z|$ is the support of $z$.

Let $Z \subset X$ be an integral closed subscheme of codimension $c$ and let $\iota: Z \to X$ be the corresponding embedding. If $Z$ itself happens to be smooth, then $\iota_* \mathcal{C}\ell_Z[Z] = \mathcal{C}\ell_X[Z]$, or equivalently

$$\mathrm{Tr}_X(\mathcal{C}\ell_X[Z] \cup u) = \mathrm{Tr}_Z(\iota^* u) \tag{7.3}$$

for all $u \in \mathrm{H}^{2(x-c)}(X_{\text{ét}}, \mathbb{Q}_\ell(x-c))$; this is [ibid., 2.3.2 and 2.3.8i]. Based on this fact, there is an entirely different construction of cycle classes. In this construction one defines weighted trace maps such that (7.3) obtains a meaning also when $Z$ is not smooth (see [ibid., 4.3]). This alternative constructions allows one to prove the following proposition. Recall that cycles $z = \sum_Z d_Z[Z]$ and $z' = \sum_Z d'_Z[Z]$ (finite sums over the closed integral subschemes) *intersect properly* if $\mathrm{Codim}(Z \cap Z') = \mathrm{Codim}(Z) + \mathrm{Codim}(Z')$ for all $Z$, $Z'$ for which $d_Z$ and $d'_{Z'}$ are both non-zero.

**7.4. Proposition.** *Let $X$, $Y$ be smooth, finite type $k$-schemes.*

*(i) The class map passes to rational equivalence.*

*(ii) If $f: X \to Y$ is flat, then $f^* \mathcal{C}\ell = \mathcal{C}\ell f^*$.*

*(iii) If $z$ and $z'$ are cycles on $X$ of codimension $c$ respectively $c'$, and if $z$ and $z'$ intersect properly, then*

$$\mathcal{C}\ell(z'') = \mathcal{C}\ell(z) \cup \mathcal{C}\ell(z')$$

*whenever $z''$ is a cycle that represents the intersection product of the cycle classes of $z$ and $z'$.*

*(iv) $\mathcal{C}\ell_{\mathrm{Spec}\,k}(\mathrm{Spec}\,k) = 1$.*

PROOF. Statement (i) is SGA $4\frac{1}{2}$-[*cycle*] Rem. 2.3.10, (iv) is by construction, (ii) is [ibid., 2.3.8ii] and (iii) is [ibid., 2.3.9]. As it may be not immediately clear that the language used in [ibid.] is compatible with ours (especially concerning intersection products), let us give the details for statement (iii).

Let $\iota: Z \hookrightarrow X$ and $\iota': Z' \hookrightarrow X$ be integral closed subschemes that are local complete intersections of codimension $c$ and $c'$, respectively. Let $\iota''$ be the immersion of $(Z \cap Z')_{\mathrm{red}}$ in $X$. Suppose $\mathrm{Codim}(Z \cap Z') = c + c'$ and let $Z_1, \ldots, Z_r$ be the irreducible components of the scheme $Z \cap Z'$. Denote by $\iota_j: Z_{j,\mathrm{red}} \hookrightarrow Z \cap Z'$ the canonical immersions. There is an open subscheme $U \subset X$ such that, for each $j$, the complement of $U \cap Z_j$ in $Z_j$ has codimension $> c + c'$, and such that $U \cap Z_{j,\mathrm{red}}$ is a local complete intersection. As it now suffices to prove the formula over $U$, we can assume $U = X$ and all $Z_{j,\mathrm{red}}$ are local complete intersections. Let $\xi_j$ be the generic point of $Z_j$ and put

$$\begin{aligned} d_j &= \sum_n (-1)^n \operatorname{length}_{\mathcal{O}_{X,\xi_j}} (h^n(\iota_* \mathcal{O}_Z \overset{\mathrm{L}}{\otimes} \iota'_* \mathcal{O}_{Z'}))_{\xi_j} \\ &= \sum_n (-1)^n \operatorname{length}_{\mathcal{O}_{X,\xi_j}} \operatorname{Tor}_{-n}^{\mathcal{O}_{X,\xi_j}} (\mathcal{O}_{Z,\xi_j}, \mathcal{O}_{Z',\xi_j}). \end{aligned}$$

Then on the one hand [ibid., 2.3.3.1, 2.3.4.1 and 2.3.8] say that $\mathcal{C}\ell[Z] \cup \mathcal{C}\ell[Z']$ is characterized by the equality

$$\mathrm{Tr}_X(\mathcal{C}\ell[Z] \cup \mathcal{C}\ell[Z'] \cup \iota''_* u) = \sum_{j=0}^{r} d_j \mathrm{Tr}_{Z_{j,\mathrm{red}}} \iota_j^*(u)$$

for all $u \in \mathrm{H}^{2\dim Z\cap Z'}(Z \cap Z', \mathbb{Q}_\ell(\dim Z \cap Z'))$. On the other hand, Serre's Tor formula ([Ful, 20.4]) says that

$$[Z] \cdot [Z'] = \sum_{j=0}^{r} d_k [Z_{j,\mathrm{red}}]$$

so that $\mathcal{C}\ell([Z] \cdot [Z'])$ is characterized by the same formula. □

**7.5. Remark.** In fact, [ibid., 2.3.8] gives a more general result than (ii). Let $f: Y \to X$ be a morphism, and let $z = \sum_i d_i [Z_i]$ be a cycle on $X$ of codimension $c$ such that $f^{-1} Z_i$ is either empty or of codimension $c$. Then $f^* \mathcal{C}\ell(z) = \mathcal{C}\ell(f^! z)$ in $\Gamma(Y, \mathcal{H}^{2c}_{f^{-1}|z|} \mathbb{Q}_\ell(c))$.

Now consider a smooth, finite type Deligne-Mumford $k$-stack $\mathcal{X}$ and an integral closed substack $\mathcal{Z}$ of codimension $c$. Choose an étale presentation $P: X \to \mathcal{X}$ by a finite type $k$-scheme and form the diagram

$$X \times_{\mathcal{X}} X \underset{p_2}{\overset{p_1}{\rightrightarrows}} X \overset{P}{\longrightarrow} \mathcal{X}.$$

As $p_1$ and $p_2$ are flat the above proposition gives

$$p_1^* \mathcal{C}\ell[P^{-1}\mathcal{Z}] = \mathcal{C}\ell[(Pp_1)^{-1}\mathcal{Z}] = \mathcal{C}\ell[(Pp_2)^{-1}\mathcal{Z}] = p_2^* \mathcal{C}\ell[P^{-1}\mathcal{Z}].$$

By the sheaf property, we can now glue to a local class $\mathcal{C}\ell(\mathcal{Z}) \in \Gamma(\mathcal{X}, \mathcal{H}^{2c}_{|\mathcal{Z}|} \mathbb{Q}_\ell(c))$. This extends by linearity to arbitrary cycles of codimension $c$. Since two rationally equivalent cycles on $\mathcal{X}$ remain equivalent when pulled-back along a chart of $\mathcal{X}$, the local class passes to rational equivalence. Corollary 7.2 then results in the cycle map on $\mathcal{X}$

$$A^c(\mathcal{X}) \overset{\mathrm{Cl}_{\mathcal{X}}}{\longrightarrow} \mathrm{H}^{2c}(\mathcal{X}_{\bar{k},\text{ét}}, \mathbb{Q}_\ell(c))^{G_k}.$$

The next proposition follows immediately from the one above, since we can localize to schemes.

**7.6. Proposition.**

*(i) The cycle maps commute with flat pull-backs.*

*(ii) If $z$ and $z'$ are cycles on $\mathcal{X}$ that intersect properly, then*

$$\mathrm{Cl}_{\mathcal{X}}([z]\cdot[z'])=\mathrm{Cl}_{\mathcal{X}}([z])\cup\mathrm{Cl}_{\mathcal{X}}([z']).\qquad\square$$

The big open problem remains to see that the cycle map commutes with products of cycles that do not intersect properly.

**7.7. Remarks.**

1. *Étale cohomology of proper, smooth schemes is a Weil cohomology.* Chow's moving lemma (see [Rob]) says that any two cycle classes on a smooth, projective scheme can be represented by cycles that intersect properly. By Proposition 7.4, the cycle map commutes with products of cycles that intersect properly. Therefore, étale cohomology of smooth, projective varieties is a Weil cohomology (cf., SGA $4\frac{1}{2}$-[*cycle*] 2.3, or section 7 of [Lau]). However, for proper, smooth schemes that are not necessarily projective the author does not know any reference for the well-known fact that étale cohomology is a Weil cohomology. The following proof was suggested by Illusie. See the second remark below for an alternative strategy.

Let $X$ be a non-singular, proper scheme over a field of characteristic different from $\ell$. Let $K^0(X)$ denote the Grothendieck group of finite rank vector bundles on $X$, tensored with $\mathbb{Q}$. It has a ring structure that is induced by the tensor product. Let $\mathrm{ch}_A: K^0(X)\to A(X)$, resp. $\mathrm{ch}_H: K^0(X)\to \mathrm{H}(X_{\text{ét}},\mathbb{Q}_\ell)$, be the Chern character to the Chow-ring, resp. to cohomology (see [Ful], §15.1; cf., [HAG], appendix A). (Note that, by our conventions, $A(X)$ is a $\mathbb{Q}$-algebra.) These maps are ring homomorphisms. By a theorem of Borelli (see [Bor]) $K^0(X)$ is isomorphic to the Grothendieck group $K(X)$ of coherent sheaves on $X$ via the natural map $K^0(X)\to K(X)$. The Chern character mapping to the Chow ring $\mathrm{ch}_A: K(X)\to A(X)$ then is an isomorphism (see [Ful], example 15.2.6b).

Now let $\mathrm{Cl}'=\mathrm{ch}_H\circ\mathrm{ch}_A^{-1}$. It is a ring homorphism $A(X)\to \mathrm{H}(X_{\text{ét}},\mathbb{Q}_\ell)$ that is functorial for pull-backs $f^*$ whenever $f: Y\to X$ is a morphism of non-singular, proper schemes. The map $\mathrm{Cl}'$ satisfies all the properties for the class map in the definition of a Weil cohomology of Definition 3.3. In fact, if $X$ is projective, then $\mathrm{Cl}=\mathrm{Cl}'$ by the universal properties of the Chern characters ([Ful] §15.1). To see this, note that we already know (by the moving lemma) that $\mathrm{Cl}\circ\mathrm{ch}_A$ is a ring homorphism; the definition of the cycle map in SGA $4\frac{1}{2}$-[*cycle*] 2.1, then shows that $\mathrm{Cl}\circ\mathrm{ch}_A$ maps a line bundle to the class in cohomology via the map obtained from the Kummer sequence; so we arive at the definition of $\mathrm{ch}_H$, showing that for line bundles $\mathrm{Cl}\circ\mathrm{ch}_A=\mathrm{ch}_H$. Using an alteration, we can reduce to the projective case.

1bis. A more elegant approach, still using line bundles, was suggested to the author by Ben Moonen. Again, let $X$ be a proper, non-singular scheme over a field of characteristic different from $\ell$. Let $Z$ be an arbitrary cycle class, but suppose $Z'$ is the chern class of a vector bundle $E$ on $X$. We will show that $\mathrm{Cl}(Z\cdot Z')=\mathrm{Cl}(Z)\cup\mathrm{Cl}(Z')$; but since the rational Chow group of $X$ is generated by chern classes of vector bundles this then settles the general case.

Let $F$ be the flag variety of the projective bundle associated to $E$ and let $p{:}F \to X$ be the canonical map. Then $p^{-1}E$ has a filtration by subbundles with line bundle quotients and $p^*$ is an injective map both on Chow groups and on cohomology—this is the *splitting principle* (see [Ful, pp. 51–54]). So $p^*Z'$ is a polynomial in divisor classes. But intersection with a divisor class commutes with the class map, as is observed in SGA $4\frac{1}{2}$[cycle] 2.1.1. Therefore,

$$\mathrm{Cl}(p^*(Z\cdot Z')) = \mathrm{Cl}(p^*Z\cdot p^*Z') = \mathrm{Cl}(p^*Z)\cup \mathrm{Cl}(p^*Z').$$

Since $p$ is flat, we may interchange Cl and $p^*$ by the above proposition. The statement now follows from the fact that $p^*$ is injective.

1ter. The above proof may be generalizable to stacks. It is true that a non-singular Deligne-Mumford stack $\mathscr{X}$ for which the map $K^0_{\mathbb{Z}}(\mathscr{X}) \to K_{\mathbb{Z}}(\mathscr{X})$ is an isomorphism is necessarily a quotient stack (this is Theorem 2.14 of [EHKV]), but as we work with rational chow groups this may not be an obstruction in the present case. (The author likes to thank Ben Moonen again for pointing this out to him.) It may be that the above approach is amendable to quotient stacks, but for such stacks the equivariant intersection theory of Edidin and Graham ([E-G]) seems the more natural and powerful approach.

2. An alternative strategy for proving that étale cohomology of a smooth and proper scheme is a Weil cohomology is to adapt the approach of Fulton ([Ful], ch. 19; cf., [Pet] for a compact exposure) which uses orientation classes of regularly imbedded subschemes. The biggest problem in adapting this strategy to stacks seems to be the definition of orientation classes of schemes that are only locally regularly imbedded (in the sense defined above), and to prove that these classes satisfy the desirable properties. Another obstruction, blowing up along a stack that is only locally embedded, seems to be overcome by Kresch in [Kr].

3. Since we are not able to prove that the cycle map commutes with arbitrary intersections, we can also not prove the equivalent statement that the cycle map commutes with arbitrary pullbacks. What we can prove is that the cycle map commutes with *push-forwards*. We will sketch the proof below. The proof relies heavily on the functorial properties of Poincaré duality, which follow from the naturality of the map (6.5). Since we have not explicitly checked this, we state the following lemma and proposition inside a remark. We will not use them in the rest of this text.

**Claim.** *Let $\mathscr{X}$ be a connected, smooth, proper Deligne-Mumford k-stack. Let $i{:}\mathscr{Z} \hookrightarrow \mathscr{X}$ be an integral closed substack of dimension z. Then*

$$\mathrm{Tr}_{\mathscr{Z}}(i^*a) = \mathrm{Tr}_{\mathscr{X}}(\mathrm{Cl}_{\mathscr{X}}[\mathscr{Z}]\cup a)$$

*for all $a \in \mathrm{H}^{2z}(\mathscr{X}_{\bar{k},\text{ét}}, \mathbb{Q}_\ell(z))$.*

Sketch of proof. The proof is essentially the one in [Laf, A.11]: reduce the problem to an étale covering by using the formalism of Poincaré duality and then employ the analogous equality for schemes given in SGA $4\frac{1}{2}$. Let us give some more details.

Let $x = \dim(\mathscr{X})$ and $c = x - z$. Consider $\mathscr{C}\ell(\mathscr{Z}) \in \Gamma(\mathscr{X}, \mathscr{H}^{2c}_{|\mathscr{Z}|}\mathbb{Q}_\ell(c))$, the local class of $\mathscr{Z}$ defined above. The cup product with $\mathscr{C}\ell(\mathscr{Z})$ (see SGA $4\frac{1}{2}$-[*cycle*] 1.2.2) defines a morphism

$$i_*\mathbb{Q}_\ell \longrightarrow \mathscr{H}^{2c}_{|\mathscr{Z}|}\mathbb{Q}_\ell(c).$$

Let $f\colon X \to \mathcal{X}$ be a non-empty étale chart and let $f'\colon Z = \mathcal{Z} \times_{\mathcal{X}} X \to \mathcal{Z}$ and $i'\colon Z \hookrightarrow X$ be the canonical maps. Then the pull-back of the above map to $X$

$$i'_*\mathbb{Q}_\ell \longrightarrow \mathcal{H}^{2c}_{|Z|}\mathbb{Q}_\ell(c)$$

is a cup product with $\mathcal{C}\ell_X(Z) = f^*\mathcal{C}\ell_{\mathcal{X}}(\mathcal{Z})$. It induces a map on the étale cohomology over $\overline{k}$

$$\mathrm{H}^{2z}_c(Z, \mathbb{Q}_\ell(z)) \longrightarrow \mathrm{H}^{2x}_c(X, \mathbb{Q}_\ell(x))$$

and, as shown in [ibid., 2.3.1], the composition of this map with $\mathrm{Tr}_X$ equals $\mathrm{Tr}_Z$.

The above morphism and its counterpart $\mathrm{H}^{2z}(\mathcal{Z}, \mathbb{Q}_\ell(z)) \to \mathrm{H}^{2x}(\mathcal{X}, \mathbb{Q}_\ell(x))$ for $\mathcal{X}$ sit in a square with the push forward morphisms, which commutes by the projection formula and the fact that by the above proposition $f^*\,\mathrm{Cl}_{\mathcal{X}}[\mathcal{Z}] = \mathrm{Cl}_X[Z]$. The dual of this square is

$$\begin{array}{ccc} \mathrm{H}^{2z}_c(Z, \mathbb{Q}_\ell(z))^\vee & \longleftarrow & \mathrm{H}^{2x}_c(X, \mathbb{Q}_\ell(x))^\vee \\ {\scriptstyle f'^\vee_*}\uparrow & & \uparrow{\scriptstyle f^\vee_*} \\ \mathrm{H}^{2z}(\mathcal{Z}, \mathbb{Q}_\ell(z))^\vee & \longleftarrow & \mathrm{H}^{2x}(\mathcal{X}, \mathbb{Q}_\ell(x))^\vee \end{array} .$$

Now $f^\vee_*(\mathrm{Tr}_{\mathcal{X}}) = \mathrm{Tr}_X$ and $f'^\vee_*(\mathrm{Tr}_{\mathcal{Z}}) = \mathrm{Tr}_Z$ and we have already observed that the top arrow maps $\mathrm{Tr}_X$ to $\mathrm{Tr}_Z$. Since $f'^\vee_*$ is injective (for it is just $f'^*\colon \mathrm{H}^0(\mathcal{Z}, \mathbb{Q}_\ell) \to \mathrm{H}^0(Z, \mathbb{Q}_\ell)$ by the isomorphism of Poincaré duality), the bottom arrow maps $\mathrm{Tr}_{\mathcal{X}}$ to $\mathrm{Tr}_{\mathcal{Z}}$. Composing this bottom arrow with the dual of $i^*$ gives the desired equality. □

We want to show that $\mathrm{Cl}_{\mathcal{Y}} f_* = f_*\,\mathrm{Cl}_{\mathcal{X}}$, where $f\colon \mathcal{X} \to \mathcal{Y}$ is a map between smooth, proper Deligne-Mumford $k$-stacks, which we may also assume to be connected.

Let $\mathcal{Z} \hookrightarrow \mathcal{X}$ be a closed, integral substack of dimension $z$, and let $\mathcal{Z}'$ be the image of $\mathcal{Z}$ under $f$ with the reduced closed substack structure. Let $f'\colon \mathcal{Z} \to \mathcal{Z}'$ and $i'\colon \mathcal{Z}' \hookrightarrow \mathcal{Y}$ be the canonical maps. Then in the Chow group of $\mathcal{Y}$ we have $f_*[\mathcal{Z}] = d[\mathcal{Z}']$, where $d$ is the generic degree of $\mathcal{Z}$ over $\mathcal{Z}'$ defined in [Vis, p. 620] (note that $d$ is a non-negative rational number, which equals 0 if $\dim \mathcal{Z} \neq \dim \mathcal{Z}'$). Now the formalism of Poincaré duality characterizes $f_*\,\mathrm{Cl}_{\mathcal{X}}[\mathcal{Z}]$ by the statement that

$$\mathrm{Tr}_{\mathcal{Y}}(f_*\,\mathrm{Cl}_{\mathcal{X}}[\mathcal{Z}] \cup a) = \mathrm{Tr}_{\mathcal{X}}(\mathrm{Cl}_{\mathcal{X}}[\mathcal{Z}] \cup f^*a) \qquad \text{for all } a \in \mathrm{H}^{2z}(\mathcal{X}_{\overline{k},\text{ét}}, \mathbb{Q}_\ell(z))$$

(where $z = \dim(\mathcal{Z})$). So we nee to verify that, for all $a$,

$$\mathrm{Tr}_{\mathcal{X}}(\mathrm{Cl}_{\mathcal{X}}[\mathcal{Z}] \cup f^*a) = \mathrm{Tr}_{\mathcal{Y}}(d \cdot \mathrm{Cl}_{\mathcal{Y}}[\mathcal{Z}'] \cup a).$$

By the above claim, this amounts to the equality

$$\mathrm{Tr}_{\mathcal{Z}}((f')^*b) = d \cdot \mathrm{Tr}_{\mathcal{Z}'}(b) \qquad (\text{with } b = (i')^*a).$$

Again, the proof of this fact relies on the functorial properties of the trace map constructed from (6.5).

## 8. The comparison theorem via simplicial schemes

In this section we prove the comparison theorem for Deligne-Mumford stacks, Statement 0.1. The method used is entirely different from the approach indicated in the previous sections. In fact, it is based on an extension by Kisin or Tsuji of the comparison theorem to simplicial schemes. The idea is that for a smooth, proper Deligne-Mumford stack one can construct a smooth, proper simplicial scheme whose étale and de Rham cohomologies agree with that of the stack. The idea of this construction is a close variant of a construction due to Deligne ([Del74ii, 6.2]). It relies heavily on De Jong's alteration theorem [dJ], and in fact a large part of this section is concerned with working out (for stacks) a construction mentioned in the third paragraph of the introduction of De Jong's paper [ibid.].

We need the theory of hypercoverings, which provides a powerful extension of the machinery of Čech cohomology. A quick introduction to hypercoverings is given in Deligne's article about Hodge theory on singular varieties [Del74ii, §5]. An extensive treatment can be found in SGA 4 $\mathrm{V}^{\mathrm{bis}}$. Besides these classical references, there is a detailed account [Co2] by Conrad, which clarifies a lot of points which are a bit obscure in SGA. Unfortunately, none of these texts is readily applicable to stacks. All results in these references are stated in terms of topological spaces, schemes, or at best (ordinary) categories that share some good properties with the category of schemes, like having fibre products.

Probably the best thing would be to rewrite SGA 4 $\mathrm{V}^{\mathrm{bis}}$ completely in the language of 2-categories. This can no doubt be done in a straightforward, though tedious, way and all texts known to the author that deal with hypercoverings of stacks make the implicit assumption that this is indeed possible. We will not carry out this generalisation here, nor use it, although we sometimes refer to it in the remarks. Instead, below we will describe an ad hoc method that gives a shortcut to the desired comparison result. For this we introduce the 'fibre category of schemes over a stack $\mathscr{X}$'. This category has the desired fibre products, but as it lacks a final object (if $\mathscr{X}$ is not representable by a scheme) we have to be careful when applying results from the literature.

### *Simplicial objects*

We start by recalling from SGA 4, $\mathrm{V}^{\mathrm{bis}}$ and V 7, the general machinery of simplicial objects. Let $\Delta$ denote the category whose objects are the sets $[k] := \{0, \dots, k\}$ with $k \geq 0$ and whose morphism are the maps preserving the relation $\leq$. For $n \in \mathbb{Z}$ denote by $\Delta_n$ the full subcategory with objects $[k]$ with $0 \leq k \leq n$.

Let $B$ be a category. A *simplicial object of* $B$ is a functor $X_\bullet : \Delta^{\mathrm{op}} \to B$. The simplicial objects of $B$ form a category $\Delta B$ with natural transformations as morphisms. The following notations are standard in the literature: A simplicial object is denoted by $X_\bullet$ (or $Y_\bullet$, etc.), the image of $[k]$ by $X_k$. Let $0 \leq i \leq k$. Consider the unique surjective map $[k+1] \to [k]$ in $\Delta$ that maps $i$ and $i+1$ to the same element of $[k]$; its image under $X_\bullet$ is denoted by $s_i^k : X_k \to X_{k+1}$. Likewise, if $k \neq 0$ then $d_i^k : X_k \to X_{k-1}$ is the image of the unique injection $[k-1] \to [k]$ in $\Delta$ that misses $i \in [k]$. A simplicial scheme $X_\bullet$ is completely determined by the various $X_i$, $s_i^k$ and $d_i^k$.

An $n$-truncated simplicial object of $B$ is a functor $\Delta_n^{\mathrm{op}} \to B$. We use similar notations as for simplicial objects and identify $\Delta_0 B$ with $B$. The inclusion functor $\Delta_n \hookrightarrow \Delta$ defines a functor $\mathrm{sk}_n : \Delta B \to \Delta_n B$, called the *n-skeleton*. A right adjoint is denoted by $\mathrm{cosk}_n^B$ and called the *coskeleton*. We drop the label $B$ or $n$ if no confusion seems likely. It exists if $B$ admits non-empty finite projective limits (for an explicit formula, see SGA 4 $\mathrm{V}^{\mathrm{bis}}$ 3.0.1.2). For example, $(\mathrm{cosk}_0^B X)_i$ is the product of $i+1$ copies of $X$.

**8.1. Lemma.** *Let $n \geq 1$. If $B$ admits finite (non-empty) fibre products and equalizers, then $\mathrm{cosk}_n^B$ exists. Moreover, the formation of the n-coskeleton commutes (up to canonical isomorphism) with any functor $B \to B'$ that commutes with finite fibre products and equalizers.*

Proof. Let $X_\bullet$ be an $n$-truncated simplicial object of $B$. First we construct the objects $(\mathrm{cosk}_n^B X_\bullet)_p$. Since $\mathrm{sk}_n \, \mathrm{cosk}_n^B X_\bullet$ is isomorphic to $X_\bullet$, we may suppose $p \geq n$. From SGA 4 $\mathrm{V}^{\mathrm{bis}}$ 3.0.1 we copy the following description of the $n$-coskeleton. Let $p \geq 0$ be an integer. Let $\Delta^+$ be the category with the same objects as $\Delta$, but whose maps are the ones that preserve the relation $<$ (not just $\leq$). Let $\Delta^+_{[p]}$ be the fibre category over $[p]$ and define $\Delta^+_{n[p]}$ to be its full subcategory whose objects $[q] \to [p]$ satisfy $0 \leq q \leq n$. Then $(\mathrm{cosk}_n^B X_\bullet)_p$ is the projective limit of the composition $\Delta^{+\mathrm{op}}_{n[p]} \to \Delta^{\mathrm{op}} \to B$, where the first map is the forgetful map and the second is given by $X_\bullet$. Now $\Delta^+_{n[p]}$, with $n > 0$, has the following three properties: (i) it is connected; (ii) given two objects $\lambda : [q] \to [p]$ and $\lambda' : [q'] \to [p]$ with $q \geq q'$, the only possible map $\lambda \to \lambda'$ is the identity; (iii) for every object $[q] \to [p]$ there exists an arrow to some $[n] \to [p]$. From this description, it follows that the projective limit can be calculated entirely in terms of fibred products and equalizers.

It remains to check functoriality. Let $\Delta_{n[p]}$ be the full subcategory of $\Delta_{[p]}$ whose objects $[q] \to [p]$ satisfy $0 \leq q \leq n$. There is a canonical isomorphism ([ibid.]):

$$\varprojlim_{\Delta_{n[p]}} X_\bullet \xrightarrow{\sim} \varprojlim_{\Delta^+_{n[p]}} X_\bullet.$$

If $\lambda : [p] \to [p']$ is a morphism in $\Delta$, it induces a functor $\Delta_{n[p]} \to \Delta_{n[p']}$, which induces the desired map

$$(\mathrm{cosk}_n \mathcal{X}_\bullet)_{p'} = \varprojlim_{\Delta_{n[p']}} \mathcal{X}_\bullet \longrightarrow \varprojlim_{\Delta_{n[p]}} \mathcal{X}_\bullet = (\mathrm{cosk}_n \mathcal{X}_\bullet)_p,$$

using the universal property of projective limits. □

**8.2. Remark.** Let us give the standard example of a simplicial object in the theory of algebraic stacks (see [L-MB, 12.4]). Let $\mathscr{X}$ be an algebraic stack and consider a presentation $P: X \to \mathscr{X}$ by a scheme $X$. Then one forms a simplicial object in the category of schemes

$$\cdots X \times_{\mathscr{X}} X \times_{\mathscr{X}} X \quad \begin{matrix} \xrightarrow{d_0^2} \\ \xleftarrow{s_1^1} \\ \xrightarrow{d_1^2} \\ \xleftarrow{s_0^1} \\ \xrightarrow{d_2^2} \end{matrix} \quad X \times_{\mathscr{X}} X \quad \begin{matrix} \xrightarrow{d_0^1} \\ \xleftarrow{s_0^0} \\ \xrightarrow{d_1^1} \end{matrix} \quad X,$$

where $d_i^k$ is the projection that forgets the $(i+1)$th coordinate, while $s_i^k$ is the map that is the diagonal from the $(i+1)$th coordinate to coordinates $i+1$ and $i+2$, and which is the identity elsewhere. In the obvious generalisation of simplicial objects to 2-categories, this simplicial object would be the '0th coskeleton of $X$ relative to $\mathscr{X}$' and a 'hypercovering of $\mathscr{X}$'.

*Simplicial schemes*

Let $L$ be a field. We will denote by $(\mathrm{Sch}/L)$ the category of schemes of finite type over $L$. A *simplicial L-scheme* is a simplicial object of $(\mathrm{Sch}/L)$. If $X_\bullet$ is a simplicial $L$-scheme and $L \subset L'$ is a field extension, one forms the simplicial $L'$-scheme $(X_\bullet)_{L'}$ in the obvious way.

Let $X_\bullet$ be a simplicial $L$-scheme. One can construct a topos denoted $(X_\bullet)_{\text{ét}}$. Let us first loosely describe what it looks like. An object of this topos is a collection of sheaves $(\mathscr{F}^k \in \mathrm{ob}(X_k)_{\text{ét}})_{k \geq 0}$ and transition morphisms $\mathscr{F}^k \to X_\bullet(i)_* \mathscr{F}^{k'}$, for each map $i: [k] \to [k']$ in $\Delta$, that are compatible with the composition in $\Delta$. By abuse of language, we call these objects *sheaves on* $X_\bullet$.

For a precise definition of $(X_\bullet)_{\text{ét}}$, we will refer to SGA 4 $\mathrm{V}^{\mathrm{bis}}$. The reader who wishes so can safely skip this paragraph. What we call $(X_\bullet)_{\text{ét}}$ is in [ibid.] designated by $\underline{\Gamma}(\overline{X_\bullet})$, a notation which brings together a manifold of constructions. First we need an auxiliary category $\mathscr{E}$ whose objects are pairs $(X, \mathscr{F})$ consisting of a finite type $L$-scheme $X$ and an object $\mathscr{F}$ of $X_{\text{ét}}$, and whose morphisms $(X, \mathscr{F}) \to (Y, \mathscr{G})$ are pairs consisting of a morphism of schemes $f: X \to Y$ and a map $\mathscr{G} \to f_* \mathscr{F}$ in $Y_{\text{ét}}$. There is a canonical functor $\mathscr{E} \to (\mathrm{Sch}/L)$, whose fibres are toposes. Furthermore, any morphism in $(\mathrm{Sch}/L)$ defines pull-back and push-forward morphisms of toposes in a natural way. This means that $\mathscr{E}$ is a category bifibred in duals of toposes ([ibid., Déf. 1.2.2]). Then $\overline{X_\bullet} = (\Delta^{\mathrm{op}} \times_{(\mathrm{Sch}/L)} \mathscr{E})^{\mathrm{op}}$ is a category bifibred over $\Delta$ and $\underline{\Gamma}(\overline{X_\bullet})$ is defined in [ibid., 1.2.8] to be the functor category $\mathrm{Hom}_\Delta(\Delta, \overline{X_\bullet})$. A straightforward check shows that this object satisfies the description of $(\mathscr{X}_\bullet)_{\text{ét}}$ given above. The important point is that it is a topos, a fact that is proved in [ibid., 1.2.12].

Having a topos, we can form the derived categories $\mathrm{D}^+((\mathscr{X}_\bullet)_{\text{ét}}, \mathbb{Z}/(\ell^n))$ of bounded below complexes of $\mathbb{Z}/(\ell^n)$-modules for $n \geq 1$. Then $\mathrm{D}_c^+((X_\bullet)_{\text{ét}}, \mathbb{Z}/(\ell^n))$ denotes the subcategory of objects $x$ whose homologies $H^n(x)$ (for $n \in \mathbb{Z}$) form a constructible module on each $X_i$ (cf. SGA 4 $\mathrm{V}^{\mathrm{bis}}$ 2.4.0–2.4.2, where another auxiliary category $G$ is introduced for this purpose, which in our case must be defined as the subcategory of constructible modules). One can now form the derived category $\mathrm{D}_c^+((X_\bullet)_{\text{ét}}, \mathbb{Q}_\ell)$ of constructible $\ell$-adic sheaves in the usual way (so projective systems of objects of $\mathrm{D}_c^+((X_\bullet)_{\text{ét}}, \mathbb{Z}/(\ell^n))$).

In an analogous way, if $L = \mathbb{C}$ we can form the derived category of sheaves of abelian groups on $X_\bullet$, using the analytic topology.

We can now form the cohomology of a constructible $\ell$-adic sheaf $\mathscr{F}^\bullet$ on $X_\bullet$. We will need the particular spectral sequence (see [Del74ii, 5.2.3.2])

$$E_1^{p,q} = \mathrm{H}^q((X_p)_{\overline{K},\text{ét}}, \mathbb{Q}_\ell) \qquad \Rightarrow \qquad \mathrm{H}^{p+q}((X_\bullet)_{\overline{K},\text{ét}}, \mathbb{Q}_\ell)$$

and the analogous one for the analytic topology for simplicial schemes over $\mathbb{C}$.

**8.3. Lemma.** *Let $f_\bullet: X_\bullet \to Y_\bullet$ be a morphism of simplicial schemes. Assume that*

*(i) $f_0: X_0 \to Y_0$ is an isomorphism,*

*(ii) $f_1: X_1 \to Y_1$ is proper and surjective,*

*(iii) the natural map $Y_\bullet \to \operatorname{cosk}_1 \operatorname{sk}_1 Y_\bullet$ is an isomorphism,*

*(iv) for all $n > 0$, the natural map $X_{n+1} \to (\operatorname{cosk}_n \operatorname{sk}_n X_\bullet)_{n+1}$ is proper and surjective.*

*Then the canonical map*

$$\mathrm{H}^\bullet((\mathscr{Y}_\bullet)_{\text{ét}}, \mathbb{Q}_\ell) \xrightarrow{\sim} \mathrm{H}^\bullet((\mathscr{X}_\bullet)_{\text{ét}}, \mathbb{Q}_\ell)$$

*is an isomorphism.*

PROOF. This result follows from a theorem by Conrad, [Co2, Thm 7.22], since proper and surjective morphisms between schemes are of universal cohomological descent relative to constructible $\ell$-adic sheaves. □

*Hypercoverings*

Let $\mathscr{X}$ be a Deligne-Mumford stack of finite type over $L$. Denote by $(\mathrm{Sch}/\mathscr{X})$ the following category. Objects are pairs $(Y, f)$, with $Y$ a finite type $L$-scheme and with $f: Y \to \mathscr{X}$ a morphism of stacks (which is automatically representable). A morphism $(Y, f) \to (Y', f')$ is a pair $(G, \varphi)$ consisting of a morphism of schemes $G: Y \to Y'$ and a 2-morphism $\varphi: f \to f'G$. Composition of morphisms is defined in the obvious way. There is a forgetful functor $F: (\mathrm{Sch}/\mathscr{X}) \to (\mathrm{Sch}/L)$ (which in general is not faithful).

**8.4. Lemma.**

*(i) Let* $(Y, f)$ *and* $(Z, g)$ *be two objects of* $(\mathrm{Sch}/\mathcal{X})$*. Consider the 2-fibre square*

$$\begin{array}{ccc} Y \times_{\mathcal{X}} Z & \xrightarrow{\pi_2} & Z \\ {\scriptstyle \pi_1}\downarrow & \overset{\epsilon}{\nearrow\!\!\!\nearrow} & \downarrow{\scriptstyle g} \\ Y & \xrightarrow[f]{} & \mathcal{X} \end{array}$$

*in the category of stacks. Then* $(Y \times_{\mathcal{X}} Z, f\pi_1)$*, together with the maps* $(\pi_1, \mathrm{id}_{f\pi_1})$ *and* $(\pi_2, \epsilon)$*, represents a product of* $(Y, f)$ *and* $(Z, g)$ *in* $(\mathrm{Sch}/\mathcal{X})$*.*

*(ii) Finite fibre products in* $(\mathrm{Sch}/\mathcal{X})$ *exist and they commute with the forgetful functor* $F$*.*

*(iii) Equalizers in* $(\mathrm{Sch}/\mathcal{X})$ *exist and commute with* $F$*.*

*(iv) Finite sums exist in* $(\mathrm{Sch}/\mathcal{X})$ *and commute with* $F$*; they are disjoint and universal (SGA 4* $\mathrm{V}^{\mathrm{bis}}$ *5.1.0).*

Proof. (i) By [L-MB, 4.2], since $\mathcal{X}$ is Deligne-Mumford, the diagonal $\mathcal{X} \to \mathcal{X} \times_L \mathcal{X}$ is schematic. Taking the fibre product of this diagonal with $Y \times_L Z \to \mathcal{X} \times_L \mathcal{X}$, it follows that $Y \times_{\mathcal{X}} Z$ is a scheme. Given two morphisms

$$(Y, f) \xleftarrow{(A,\varphi)} (T, h) \xrightarrow{(B,\psi)} (Z, g)$$

in $(\mathrm{Sch}/\mathcal{X})$. There is a 2-commutative diagram

$$\begin{array}{ccc} T & \xrightarrow{B} & Z \\ {\scriptstyle A}\downarrow & \overset{\psi\varphi^{-1}}{\nearrow\!\!\!\nearrow} & \downarrow{\scriptstyle g} \\ Y & \xrightarrow[f]{} & \mathcal{X} \end{array}$$

in the 2-category of stacks. Since the only 2-isomorphisms between morphisms of schemes are identities, the universal property of 2-fibre products says there is a unique 1-morphism $C: T \to Y \times_{\mathcal{X}} Z$ such that $\pi_1 C = A$ and $\pi_2 C = B$ and for which $\epsilon * \mathrm{id}_C = \psi\varphi^{-1}$. So we obtain a factorisation

$$\begin{array}{ccc} & & (Y, f) \\ & \overset{(A,\varphi)}{\nearrow} & \uparrow{\scriptstyle (\pi_1, \mathrm{id}_{f\pi_1})} \\ (T, h) & \xrightarrow{(C,\varphi)} & (Y \times_X Z, f\pi_1) \\ & \underset{(B,\psi)}{\searrow} & \downarrow{\scriptstyle (\pi_2, \epsilon)} \\ & & (Z, g) \end{array}$$

and the morphism $(C,\varphi)$ is the unique one making this diagram commute.

(ii) Consider a commutative square

$$\begin{array}{ccc} (T,h) & \xrightarrow{(B,\varphi)} & (Z,g) \\ {\scriptstyle (A,\psi)}\big\downarrow & & \big\downarrow{\scriptstyle (p_2,\chi_2)} \\ (Y,f) & \xrightarrow[(p_1,\chi_1)]{} & (X,e) \end{array}$$

in $(\mathrm{Sch}/\mathscr{X})$. Applying the forgetful functor, we obtain a unique morphism $C\colon T \to Y \times_X Z$ for which $\pi_1 C = A$ and $\pi_2 C = B$. If we want to extend $C$ to a morphism $(T,h) \to (Y \times_X Z, f\pi_1)$ in $(\mathrm{Sch}/\mathscr{X})$, there is only one choice: $(C,\psi)$. Finally, one checks that

$$(\pi_2, (\mathrm{id}_{\pi_2} * \chi_2) \circ (\mathrm{id}_{\pi_1} * \chi_1)) \circ (C,\psi) = (B,\varphi),$$

as is required.

(iii) is proved in the same way as (ii).

(iv) follows immediately from the construction of sums in $(\mathrm{Sch}/L)$. □

If $X_\bullet$ is a simplicial object of $(\mathrm{Sch}/\mathscr{X})$, denote by $\widetilde{X}_\bullet = (\Delta F)X_\bullet$ the underlying simplicial scheme. Likewise, if $X_\bullet$ is an $n$-truncated simplicial object of $(\mathrm{Sch}/\mathscr{X})$, then $\widetilde{X}_\bullet$ denote the underlying $n$-truncated simplicial scheme. If we combine the previous lemma with Lemma 8.1, we obtain the next result.

**8.5. Lemma.** *Let $n \geq 1$ be an integer and let $X_\bullet$ be an $n$-truncated simplicial object of* $(\mathrm{Sch}/\mathscr{X})$. *Then* $\mathrm{cosk}_n^{(\mathrm{Sch}/\mathscr{X})} X_\bullet \simeq \mathrm{cosk}_n^{(\mathrm{Sch}/L)} \widetilde{X}_\bullet$. □

**8.6. Remark.** Had we developed a notion of simplicial object in a 2-category, then this lemma would follow from the well-known base change formula ([Del74ii, 5.1.1])

$$\mathrm{cosk}_n^{\mathscr{X}} X_\bullet = (\mathrm{cosk}_n^{L} X_\bullet) \times_{\mathrm{cosk}_n C_{\mathscr{X}}} \mathscr{X}.$$

Let $X_\bullet$ be a simplicial object of $(\mathrm{Sch}/\mathscr{X})$. If

(i) in $X_0 = (\widetilde{X}_0, f)$, the map $f\colon X_0 \to \mathscr{X}$ is surjective and proper,

(ii) for each $i \geq 0$ the natural map

$$\widetilde{X}_{i+1} \longrightarrow (\mathrm{cosk}_i^{(\mathrm{Sch}/\mathscr{X})} \mathrm{sk}_i X_\bullet)^{\sim}_{i+1}. \tag{8.7}$$

is surjective and proper,

then $X_\bullet$ is called a *hypercovering of* $\mathscr{X}$.

Recall form §2 the construction of the category $(\mathcal{DM}/L)^1$, obtained from the 2-category of Deligne-Mumford stacks by identifying all 2-isomorphic morphisms. Let $C_{\mathcal{X}}$ be the constant simplicial object of $(\mathcal{DM}/L)^1$ that is equal to $\mathcal{X}$ in all degrees, and all whose maps are the identities. To a simplicial object $X_\bullet$ of $(\mathrm{Sch}/\mathcal{X})$ one associates in a natural way a morphism in $\Delta(\mathcal{DM}/L)^1$

$$\widetilde{X}_\bullet \xrightarrow{\theta} C_{\mathcal{X}},$$

where we view $\widetilde{X}_\bullet$ as a simplicial object of $(\mathcal{DM}/L)^1$. (If $\mathcal{X}$ were a scheme, such a map $\theta$ would in turn define a simplicial object of $(\mathrm{Sch}/\mathcal{X})$.) As noted in §2, a morphism $f\colon \mathcal{Y} \to \mathcal{X}$ in $(\mathcal{DM}/L)^1$ induces a morphism of toposes $\mathcal{Y}_{\text{ét}} \to \mathcal{X}_{\text{ét}}$. Therefore $\theta$ induces a map

$$\mathrm{Mod}(\mathcal{X}_{\text{ét}}, \mathbb{Q}_\ell) \xrightarrow{\overline{\theta}^*} \mathrm{Mod}((\widetilde{X}_\bullet)_{\text{ét}}, \mathbb{Q}_\ell);$$

it sends a module $\mathcal{F}$ of $\mathcal{X}_{\text{ét}}$ to $(\theta_i^*\mathcal{F})_{i\geq 0}$, where $\theta_i$ is the map $X_i \to \mathcal{X}$ that is part of the data forming $\theta$. (We use the notation of SGA 4 $\mathrm{V}^{\mathrm{bis}}$ 2.2 here.) The map $\overline{\theta}^*$ is exact and has a right adjoint $\overline{\theta}_*$, which maps an object $(\mathcal{F}^i)_{i\geq 0}$ to the projective limit of the system of modules $(\theta_{i*}\mathcal{F}^i)_{i\geq 0}$ on $\mathcal{X}$ (this is explained in [ibid.]).

**8.8. Lemma.** *Suppose $X_\bullet$ is a hypercovering of $\mathcal{X}$. Then there is a natural isomorphism*

$$\mathrm{H}^m(\mathcal{X}_{\text{ét}}, \mathbb{Q}_\ell) \xrightarrow{\sim} \mathrm{H}^m((\widetilde{X}_\bullet)_{\text{ét}}, \mathbb{Q}_\ell).$$

PROOF. Let us first assume that $X_\bullet = \mathrm{cosk}_0^{(\mathrm{Sch}/\mathcal{X})}\mathrm{sk}_0 X_\bullet$, so that $\widetilde{X}_i$ equals the fibre product over $\mathcal{X}$ of $i+1$ copies of $\widetilde{X}_0$. There is a natural map

$$\mathbb{Q}_\ell \longrightarrow \mathrm{R}\overline{\theta}_*\mathbb{Q}_\ell$$

and it suffices to check that this is an isomorphism, or, which amounts to the same thing, that $\mathbb{Q}_\ell \simeq \mathrm{R}^0\overline{\theta}_*\mathbb{Q}_\ell$ and $\mathrm{R}^n\overline{\theta}_*\mathbb{Q}_\ell = 0$ for $n \neq 0$. But this is a local question on $\mathcal{X}$, and therefore we may pull-back along an étale morphism $g\colon U \to \mathcal{X}$ with $U$ a scheme. From SGA 4 $\mathrm{V}^{\mathrm{bis}}$ 2.5.4, combined with exactness of $g^*$, it follows that there is a spectral sequence

$$E_1^{p,q} = g^*\mathrm{R}^q(\theta_p)_*\mathbb{Q}_\ell \qquad \Rightarrow \qquad g^*\mathrm{R}^{p+q}\overline{\theta}_*\mathbb{Q}_\ell.$$

From base change results for stacks (see [L-MB, 18.5.1], although a much weaker result than this suffices) it follows that this spectral sequence coincides with that of

$$(\widetilde{X}_\bullet \times_{\mathcal{X}} U) \longrightarrow C_U.$$

Since properness and surjectiveness are preserved under base change, this is a hypercovering of *schemes*, so we are done by SGA 4 $\mathrm{V}^{\mathrm{bis}}$ 3.3.3 and 4.3.2.

Now we prove the general case. Put $Y_\bullet = (\operatorname{cosk}_0^{(\mathrm{Sch}/\mathscr{X})} \operatorname{sk}_0 X_\bullet)^\sim$ and let $f_\bullet$ be the natural map $\widetilde{X}_\bullet \to Y_\bullet$. Note that Lemma 8.5 gives

$$\operatorname{cosk}_1^{(\mathrm{Sch}/L)} \operatorname{sk}_1 Y_\bullet = (\operatorname{cosk}_1^{(\mathrm{Sch}/\mathscr{X})} \operatorname{sk}_1 \operatorname{cosk}_0^{(\mathrm{Sch}/\mathscr{X})} \operatorname{sk}_0 X_\bullet)^\sim = (\operatorname{cosk}_0^{\mathscr{X}} \operatorname{sk}_0 X_\bullet)^\sim = Y_\bullet$$

and for $n > 0$

$$\operatorname{cosk}_n^{(\mathrm{Sch}/L)} \operatorname{sk}_n \widetilde{X}_\bullet = (\operatorname{cosk}_n^{(\mathrm{Sch}/\mathscr{X})} \operatorname{sk}_n X_\bullet)^\sim.$$

Thus, we are done by Lemma 8.3. □

*Main example*

The object of the present section is to construct a hypercovering of a proper Deligne-Mumford $K$-stack by schemes that have Zariski locally a semi-stable model with respect to some finite extension of $K$. This construction is (at least in the context of schemes) well-known. By lack of a suitable reference, we will supply the details.

We adopt the definition of semi-stabilty of [dJ, 2.16]: a variety defined over a discrete valuation ring with perfect residue field is semi-stable if it is regular, if its generic fibre is smooth and the special fibre is a reduced divisor with normal crossings. Note that semi-stability is *not* stable under base change. Heavy use is made of the alteration theorem of De Jong, formulated in the next lemma (which is not its most general formulation).

**8.9. Lemma** (De Jong [dJ, 6.5]). *Let $\widetilde{K}$ be a complete discrete valuation field. Let $X$ be an integral, proper and flat scheme over $\mathcal{O}_{\widetilde{K}}$. There exists a finite extension $K'$ of $\widetilde{K}$, an integral, proper and flat $\mathcal{O}_{K'}$-scheme $X'$, and a dominant, proper and generically finite $\mathcal{O}_{\widetilde{K}}$-morphism $\varphi: X' \to X$ such that $X'$ is semi-stable.* □

Let $\mathscr{X}$ be a proper Deligne-Mumford $K$-stack. By [L-MB, 16.6.1] ('Chow's lemma for Deligne-Mumford stacks') there exists a projective $K$-scheme $X_o$ together with a proper, surjective morphism $X_o \to \mathscr{X}$. Fix a closed embedding of $X_o$ into projective $N$-space $\mathbb{P}_K^N$, for some $N \in \mathbb{Z}$, and define $X''$ as the closure of $X_o$ in $\mathbb{P}_{\mathcal{O}_K}^N$. (Alternatively, let $X_o$ be the scheme given by Theorem 0.2, consider it as a scheme over $\mathcal{O}_K$ and apply Nagata's theorem explained in [Co1] to define $X''$.) We may assume that $X''$ is flat over $\mathcal{O}_K$. Let $X$ be the disjoint union of the irreducible components of $X''$, each being given the reduced subscheme structure. Then we can apply the above lemma componentwise to obtain a finite extension $K^{(0)}$ of $K$ (which we may suppose to be contained in $\overline{K}$) and a proper and flat $\mathcal{O}_{K^{(0)}}$-scheme $X^{(0)}$. Zariski locally, this scheme is the base change of a semi-stable scheme defined over the ring of integers of some intermediate field of $K$ and $K^{(0)}$. By construction, there is a proper, surjective morphism from $X^{(0)}_{K^{(0)}}$ to $\mathscr{X}_{K^{(0)}}$.

In the construction that follows, we are going to extend $X^{(0)}$ to an $n$-truncated simplicial scheme. More precisely, we will, for each $n \geq 0$, define a finite extension $K^{(n)}$ of $K^{(0)}$ inside $\overline{K}$ and an $n$-truncated simplicial $\mathcal{O}_{K^{(n)}}$-scheme $X_\bullet^{(n)}$. These data will satisfy the following properties for all $0 \leq i \leq n$:

(i) The scheme $X_i^{(n)}$ is proper and flat over $\mathcal{O}_{K^{(n)}}$.

(ii) There exists, for each $0 \leq j \leq n$, a subscheme $NX_j$ of $X_j^{(n)}$ such that for every $0 \leq i \leq n$ the map
$$\bigsqcup NX_j^{(n)} \longrightarrow X_i^{(n)}$$
is an isomorphism; here the disjoint union is indexed by the order preserving surjections $[i] \rightarrow [j]$ in $\Delta$ with $j \leq i$ and the morphism is the one induced by these surjections. Furthermore, each $NX_j$ can, Zariski locally, be obtained as the base change of a semi-stable scheme defined over the ring of integers of an intermediate field of $K$ and $K^{(n)}$.

(iii) There exists an $n$-truncated simplicial object $Y_\bullet^{(n)}$ of $(\mathrm{Sch}/\mathcal{X}_{K^{(n)}})$ whose underlying $n$-truncated simplicial scheme $\widetilde{Y}_\bullet^{(n)}$ is isomorphic with $(X_\bullet^{(n)})_{K^{(n)}}$, and whose coskeleton is a hypercovering of $\mathcal{X}_{K^{(n)}}$.

(Statement (ii) is a bit technical. It is necessary to be able to employ the method of SGA [$\mathrm{V}^{\mathrm{bis}}$ 5.1] for defining a simplicial scheme inductively. Anyhow, it follows from this statement that each $X_i^{(n)}$ comes Zariski locally from a semi-stable scheme, which is the property we will need in the next section.)

We have already defined $K^{(0)} = K'$ and we identify $X^{(0)}$ with a 0-truncated simplicial scheme $X_\bullet^{(0)}$. Putting $Y_\bullet^{(0)}$ equal to $X_{K^{(0)}}^{(0)}$, the above properties (i)–(iii) are satisfied. Let $n > 0$ be an integer and suppose, inductively, that $K^{(i)}$ and $X_\bullet^{(i)}$ have been defined for all $0 \leq i < n$.

Let $Y_\bullet^{(n-1)}$ be the $(n-1)$-truncated simplicial object of $(\mathrm{Sch}/\mathcal{X}_{K^{(n-1)}})$ of property (iii) above. Let
$$X' = (\mathrm{cosk}^{(\mathrm{Sch}/\mathcal{O}_{K^{(n-1)}})} X_\bullet^{(n-1)})_n \qquad \text{and} \qquad Y' = (\mathrm{cosk}^{(\mathrm{Sch}/\mathcal{X}_{K^{(n-1)}})} \widetilde{Y}_\bullet^{(n-1)})_n \,.$$

There is a canonical map $Y' \rightarrow X'$. In fact, if $n > 1$ then the coskeletons are calculated entirely in terms of fibre products and equalizers (Lemma 8.1) and these agree in $(\mathrm{Sch}/K^{(n-1)})$ and $(\mathrm{Sch}/\mathcal{X}_{K^{(n-1)}})$ (Lemma 8.4), so the map $Y' \rightarrow (X')_{K^{(n-1)}}$ is the identity. If $n = 1$, then $Y' \rightarrow X'$ is the map
$$\widetilde{Y}_0^{(0)} \times_{\mathcal{X}_{K^{(n-1)}}} \widetilde{Y}_0^{(0)} \longrightarrow X_0^{(0)} \times_{\mathcal{O}_{K^{(n-1)}}} X_0^{(0)}$$
(the left hand side is obtained from Lemma 8.4 again). As in both cases the map is separated and quasi-finite, Zariski's main theorem (EGA IV 8.12.16) tells us that it factors through an open immersion with dense image $Y' \to X''$ and a proper morphism $X'' \to X'$. Like before, we let $X$ be the disjoint union of the irreducible components of $X''$, equiped with the reduced subscheme structure, and let $Y$ be the pull-back of $Y'$. Since both $X$ and $Y$ are proper, $Y = X_{K^{(n-1)}}$. Apply the above lemma of De Jong to obtain a finite extension $K^{(n)}$ of $K^{(n-1)}$, a proper, flat scheme $NX_n^{(n)}$ over $\mathcal{O}_{K^{(n)}}$, and a proper, dominant morphism $NX^{(n)} \to X'$.

For $0 \leq i < n$, define $X_i^{(n)} = (X_i^{(n-1)})_{\mathcal{O}_{K^{(n)}}}$ and $NX_i^{(n)} = (NX_i^{(n-1)})_{\mathcal{O}_{K^{(n)}}}$. Put

$$X_n^{(n)} = \bigsqcup NX_j^{(n)},$$

the disjoint union being indexed by all order preserving, surjective maps $[n] \to [j]$ with $0 \leq j \leq n$. Using property (ii) above and the canonical map $X_n^{(n)} \to X'_{\mathcal{O}_{K^{(n)}}}$, it is showed in SGA [$\mathrm{V}^{\mathrm{bis}}$ 5.1.3] that the various schemes $X_i^{(n)}$ (with $0 \leq i \leq n$) form the objects of an $n$-truncated simplicial scheme $X_\bullet^{(n)}$ that satisfies $\mathrm{sk}_{n-1} X_\bullet^{(n)} = (X_\bullet^{(n-1)})_{\mathcal{O}_{K^{(n)}}}$. The conditions (i), (ii) are satisfied by definition.

For condition (iii), define $Y'_\bullet = (X_\bullet^{(n)})_{K^{(n)}}$. By the induction hypotheses and the construction of $X_n^{(n)}$ above, there is a map

$$Y'_\bullet \longrightarrow \mathrm{sk}_n(\mathrm{cosk}_{n-1}^{(\mathrm{Sch}/\mathcal{X}_{K^{(n)}})}(Y_\bullet^{(n-1)})_{\mathcal{X}_{K^{(n)}}})^{\sim}.$$

This implies that there exists an $n$-truncated simplicial object $Y_\bullet^{(n)}$ of $(\mathrm{Sch}/\mathcal{X}_{K^{(n)}})$ whose underlying $n$-truncated simplicial scheme equals $Y'_\bullet$. As the map above is proper and surjective, $\mathrm{cosk}\, Y_\bullet^{(n)}$ is a hypercovering of $\mathcal{X}_{K^{(n)}}$. This ends the construction of $X_\bullet^{(n)}$.

The above construction is of course highly non-canonical. By going through the construction, it is a straightforward, though tedious matter to obtain the following result. Let $(K'_1, X_{\bullet 1})$ and $(K'_2, X_{\bullet 2})$ be data both satisfying the above conditions for a particular $n$ (which we suppress from the notation). There then exists a datum $(K', X_\bullet)$ such that $K'$ contains both $K'_1$ and $K'_2$, and such that there are morphisms $X_{\bullet 1} \leftarrow X_\bullet \to X_{\bullet 2}$ of truncated simplicial objects.

**8.10. Remarks.** 1. A more elaborate version of the above argument gives the following result: Consider a finite type Deligne-Mumford $K$-stack $\mathcal{X}$ that is not necessarily proper. Then for each $n \geq 0$, there exists a finite extension $K'$ of $K$, an $n$-truncated simplicial scheme $\overline{X}_\bullet$ that satisfies conditions (i) and (ii) above, and an $n$-truncated simplicial scheme $X_\bullet$ that is a hypercovering of $\mathcal{X}_{K'}$ and that, on each level, is an open subscheme of $\overline{X}_\bullet$ whose complement is a normal crossing divisor. To get this, the above construction can easily be extended, using an appropriate version of De Jong's result for the open case.

2. A simplified version of this argument gives the following. Let $\mathcal{X}$ be a Deligne-Mumford stack of finite type over a perfect field $k$. Then there exist a hypercovering $X_\bullet$ of $\mathcal{X}$, a simplicial $k$-scheme $\overline{X}_\bullet$ and an open immersion (in the obvious sense) $\widetilde{X}_\bullet \hookrightarrow \overline{X}_\bullet$, such that each $\overline{X}_i$ is projective and smooth and $\overline{X}_i \backslash \tilde{X}_i$ is a normal crossing divisor. If $k$ is a field of characteristic 0, and if one uses Hironaka's resolution of singularities instead of De Jong's alterations, this is the construction described in [Del74ii, 6.2] and SGA [$\mathrm{V}^{\mathrm{bis}}$ 5.3].

3. All constructions made can be put in a general framework, which is a straightforward extension of the one given in SGA 4 $\mathrm{V}^{\mathrm{bis}}$, end of §5.1; for this one can use the category of weakly semi-stable pairs, introduced by [Ki].

*The comparison theorem*

Although Faltings alludes to simplicial versions of his comparison theorems in the introduction to [Fa], it is only in [Ki] that one finds the generalisations of Faltings' methods to certain simplicial schemes. Tsuji [Ts98] sketches an approach to comparison theorems using simplicial methods that is totally different from Falting's method. One can find the following result in both Kisin's and Tsuji's papers (Kisin has a more general version that does not assume properness).

**8.11. Theorem** ([Ki, 2.8.2], [Ts98, 7.1.1]). *Let $K'$ be a discrete valuation field of characteristic 0 with perfect residue field of characteristic $p$, contained in $\overline{K}$. Let $\overline{X}_\bullet$ be an $n$-truncated simplicial scheme over $\mathcal{O}_{K'}$. Suppose each $\overline{X}_i$ is a proper and flat scheme over $\mathcal{O}_{K'}$ that is étale locally the base change of a semi-stable scheme coming from a subfield $K'' \subset K'$ such that $K'$ is a finite extension of $K''$. Let $X_\bullet = (\overline{X}_\bullet)_{K'}$ be the generic fibre. For each $m \in \mathbb{Z}$ there is an isomorphism of filtered $\overline{K}$-vector spaces*

$$\mathrm{H}^m((X_\bullet)_{\overline{K},\text{ét}}, \mathbb{Q}_p) \otimes_{\mathbb{Q}_p} B_{\mathrm{dR}} \xrightarrow{\sim} \mathrm{H}^m(X_\bullet, \Omega^\bullet_{X_\bullet/K}) \otimes_K B_{\mathrm{dR}}$$

*that is equivariant for the action of* $\mathrm{Gal}(\overline{K}/K)$. □

**8.12. Corollary.** *The comparison theorem (see Statement 0.1) holds for smooth, proper Deligne-Mumford $K$-stacks.*

PROOF. Let $\mathcal{X}$ be such a stack and fix an integer $m$. We need to show that there exists an isomorphism

$$\mathrm{H}^m(\mathcal{X}_{\overline{K},\text{ét}}, \mathbb{Q}_p) \otimes_{\mathbb{Q}_p} B_{\mathrm{dR}} \xrightarrow{\sim} \mathrm{H}^m(\mathcal{X}, \Omega^\bullet_{\mathcal{X}/K}) \otimes_K B_{\mathrm{dR}}$$

respecting the filtrations and the $G_K$-actions. Consider the fields $K^{(n)}$ and the truncated simplicial objects $X^{(n)}_\bullet$ and $Y^{(n)}$ defined in the previous section for the various integers $n \geq 0$. The above theorem is applicable to $X^{(n)}_\bullet$. By Lemma 8.8, there is an isomorphism of $\mathrm{Gal}(\overline{K}/K^{(n)})$-modules

$$\mathrm{H}^m(\mathrm{cosk}(\tilde{Y}^{(n)}_\bullet)_{\overline{K},\text{ét}}, \mathbb{Q}_p) \xrightarrow{\sim} \mathrm{H}^m(\mathrm{cosk}\,\mathcal{X}_{\overline{K},\text{ét}}, \mathbb{Q}_p).$$

Now there is the spectral sequence

$$E_1^{a,b} = \mathrm{H}^b(\mathrm{cosk}(\widetilde{Y}_a^{(n)})_{\overline{K},\text{ét}}, \mathbb{Q}_p) \qquad \Rightarrow \qquad \mathrm{H}^{a+b}(\mathrm{cosk}(\widetilde{Y}_\bullet^{(n)})_{\overline{K},\text{ét}}, \mathbb{Q}_p).$$

Since $E_n^{a,b} = 0$ if either $a$ or $b$ is negative, the $m$-th cohomology only depends on $E_1^{a,b}$ for $(a,b)$ in some bounded domain. It follows that $\mathrm{H}^m(\mathrm{cosk}(\widetilde{Y}_\bullet^{(n)})_{\overline{K},\text{ét}}, \mathbb{Q}_p)$ and $\mathrm{H}^m((\widetilde{Y}_\bullet^{(n)})_{\overline{K},\text{ét}}, \mathbb{Q}_p)$ are equal for $n$ sufficiently large. From now on, we fix such a large $n$.

There also is a canonical map of filtered $K^{(n)}$-vector spaces

$$\mathrm{H}^m(\mathrm{cosk}(\widetilde{Y}_\bullet^{(n)})_{\text{ét}}, \Omega^\bullet_{\mathrm{cosk}\,\widetilde{Y}_\bullet^{(n)}/K^{(n)}}) \longrightarrow \mathrm{H}^m_{\mathrm{dR}}(\mathscr{X}/K) \otimes_K K^{(n)}.$$

To see that this map is an isomorphism, we apply the following trick. We may assume $K = K^{(n)} = \mathbb{C}$ and by applying GAGA on the components of the Hodge-de Rham spectral sequence, we may replace $\widetilde{Y}_\bullet^{(n)}$ and $\mathscr{X}$ by their analytifications. Since the holomorphic de Rham complex is a resolution of the constant sheaf $\mathbb{C}$, the problem amounts in showing that the map

$$\mathrm{H}^m(\mathrm{cosk}(\widetilde{Y}_\bullet^{(n)})^{\mathrm{an}}, \mathbb{C}) \longrightarrow \mathrm{H}^m(\mathscr{X}^{\mathrm{an}}, \mathbb{C})$$

is an isomorphism. This follows from the fact that proper hypercoverings are of cohomological descent for sheaves of abelian groups on analytic stacks.

Hence we obtain an isomorphism

$$\mathrm{H}^m(\mathscr{X}_{\overline{K},\text{ét}}, \mathbb{Q}_p) \otimes_{\mathbb{Q}_p} B_{\mathrm{dR}} \xrightarrow{\sim} \mathrm{H}^m_{\mathrm{dR}}(\mathscr{X}/K) \otimes_K B_{\mathrm{dR}}$$

of filtered $\overline{K}$-vector spaces, equivariant for the action of $\mathrm{Gal}(\overline{K}/K^{(n)})$. It remains to show that the isomorphism is in fact equivariant for the action of the bigger group $\mathrm{Gal}(\overline{K}/K)$. To prove this one first notes that this isomorphism is independent of the choice of $K^{(n)}$, $X_\bullet^{(n)}$ and $Y_\bullet^{(n)}$, since any two such data can be covered by a third. Having established this, one employs the trick of [Ts02, A7]. □

## 9. Appendix: open and closed substacks

In order to treat cohomology with supports, it is useful to have functors like 'extension by zero of sheaves on an open substack' or 'subsheaves with support in a closed substack', analogous as for schemes (see for instance [Mi, pp. 73–76]). The nice thing is that such constructions and their most important properties follow from general topos theory; this is described in SGA 4 IV 9. In order to apply these general results, one needs to link the notions of open substacks, closed substacks and complements to the analogous notions in topos theory. For schemes this is done in SGA 4 VIII 6 and the purpose of this section is to extend these results to stacks.

Let $\mathcal{X}$ be a (quasi-separated) Deligne-Mumford stack defined over a quasi-separated base scheme $S$. By $\mathcal{X}_{\text{ét}}$ we denote the étale topos associated to $\mathcal{X}$ (see section 2). Recall that an open, resp. closed, substack $\mathcal{Y}$ of $\mathcal{X}$ is a sub-$S$-stack $\mathcal{Y}$ of $\mathcal{X}$ such that the canonical map $\mathcal{Y} \to \mathcal{X}$ is an open, resp. closed, immersion ([L-MB, 3.14]). Open substacks of $\mathcal{X}$ correspond canonically to open subsets of the Zariski space associated to $\mathcal{X}$, and this correspondence preserves inclusions ([ibid., 5.4 and 5.5]). Likewise, a closed substack induces a closed subset of the Zariski space, and each closed subset is induced (non-uniquely) in this way ([ibid., 5.6.1]).

Let $j\colon \mathcal{U} \to \mathcal{X}$ be a representable morphism of Deligne-Mumford stacks. For each étale open $(V, v)$, let $\widetilde{\mathcal{U}}(V, v)$ be the set of sections of the canonical projection of $\mathcal{U} \times_{\mathcal{X}} V$ to $V$. If $(\varphi, \alpha)\colon (W, w) \to (V, v)$ is a morphism in the étale site of $\mathcal{X}$, the universal property of the fibre product determines a map $\widetilde{\mathcal{U}}(V, v) \to \widetilde{\mathcal{U}}(W, w)$. Since the map $\mathcal{U} \times_{\mathcal{X}} V \to V$ is a map of algebraic spaces, it follows that $\widetilde{\mathcal{U}}$ is a sheaf on the étale site of $\mathcal{X}$. We will use this construction below if $j$ is an open immersion.

**9.1. Proposition** (compare SGA 4 VIII 6.1). *The mapping $\mathcal{U} \to \widetilde{\mathcal{U}}$ defines an inclusion preserving, bijective map from the set of Zariski opens of $\mathcal{X}$ to the set of opens of $\mathcal{X}_{\text{ét}}$, the étale topos of $\mathcal{X}$.*

PROOF. Let $\mathcal{U}$ be an open substack of $\mathcal{X}$ and let $j\colon \mathcal{U} \to \mathcal{X}$ be the corresponding open immersion. By definition $\widetilde{U}$ is open if it is a subobject of the final object in $\mathcal{X}_{\text{ét}}$. This is the case, since for any étale open $(V, v)$ the morphism $\mathcal{U} \times_{\mathcal{X}} V \to V$ is an open immersion and therefore $\widetilde{\mathcal{U}}(V, v)$ consists of at most one object.

That $\mathcal{U} \mapsto \widetilde{U}$ preserves inclusion is clear, thus it is an order preserving map between the set of Zariski opens and the set of opens in $\mathcal{X}_{\text{ét}}$ (and hence injective). To see it is bijective, let $\mathcal{F}$ be a open object of $\mathcal{X}_{\text{ét}}$. If $v\colon V \to \mathcal{X}$ is an étale morphism, then $v(|V|) \subset |\mathcal{X}|$ is Zariski open ([L-MB, 5.6]). Hence the union of subsets $v(|V|)$ over all étale opens $(V, v)$ such that $\mathcal{F}(V, v) \neq \emptyset$ is open, so it corresponds to some open substack $\mathcal{U}$ of $\mathcal{X}$. By construction $\mathcal{F} = \widetilde{\mathcal{U}}$. □

We immediately obtain (compare SGA 4 VIII 6.2) that the morphism $(j_*, j^{-1})\colon \mathcal{U}_{\text{ét}} \to \mathcal{X}_{\text{ét}}$ realizes $\mathcal{U}_{\text{ét}}$ as an open subtopos of $\mathcal{X}_{\text{ét}}$, associating $\mathcal{U}_{\text{ét}}$ to $(\mathcal{X}_{\text{ét}})_{/\widetilde{U}}$. A little more generally, an arbitrary open immersion $j\colon \mathcal{U} \to \mathcal{X}$ defines an open embedding (SGA 4 IV 9.2) $(j_*, j^{-1})\colon \mathcal{U}_{\text{ét}} \to \mathcal{X}_{\text{ét}}$. In fact, it suffices to check that $j_*$ is fully faithful, which is straightforward.

Now we turn our attention to closed substacks. First recall from [L-MB, 12.2.2] that $\mathcal{X}_{\text{ét}}$ has sufficiently many points: if $r\colon \mathcal{F} \to \mathcal{G}$ is a morphism of sheaves on the étale site of $\mathcal{X}$, then $r$ is an isomorphism if and only if for every étale open $(U,u)$ and every geometric point $\xi$ of $U$ the induced map $(\mathcal{F}_{U,u})_\xi \to (\mathcal{G}_{U,u})_\xi$ is an isomorphism.

Now let $i:\mathcal{Z} \to \mathcal{X}$ be a closed immersion. Then for any sheaf $\mathcal{F}$ on $\mathcal{Z}$ the canonical map $i^{-1}i_*\mathcal{F} \to \mathcal{F}$ is an isomorphism. Indeed, by the above it suffices to show this on points and hence it follows from the analogous results for schemes—see SGA 4 VIII 5.7. Hence $i_*$ is fully faithful, so $(i_*, i^{-1})$ is an embedding (SGA 4 IV 9.1.1) of $\mathcal{Z}_{\text{ét}}$ in $\mathcal{X}_{\text{ét}}$.

Let $\mathcal{U}$ be an open substack of $\mathcal{X}$. There exists a unique closed substack $\mathcal{Z}$ such that for each presentation $X \to \mathcal{X}$ the pull-back of $\mathcal{Z}$ to $X$ is the reduced closed algebraic subscheme that is complementary to the pull-back of $U$ ([L-MB, 4.10]). Denote by $j:\mathcal{U} \to \mathcal{X}$ and $i:\mathcal{Z} \to \mathcal{X}$ the canonical embeddings.

**9.2. Proposition** (compare SGA 4 VIII 6.3). *The closed subtopos $\mathcal{Z}_{\text{ét}}$ of $\mathcal{X}_{\text{ét}}$ is complementary to the open subtopos $\mathcal{U}_{\text{ét}}$. This means (SGA 4 IV 9.3.4) that for any sheaf $\mathcal{F}$ on $\mathcal{X}$ we have: $j^{-1}\mathcal{F}$ is a final object of $\mathcal{U}_{\text{ét}}$ if and only if $\mathcal{F} \simeq i_*\mathcal{G}$ for some sheaf $\mathcal{G}$ on $\mathcal{Z}$.*

Proof. Sufficiency follows from the fact that $\mathcal{U} \times_{\mathcal{X}} \mathcal{Z}$ is the empty stack. Suppose that $j^{-1}\mathcal{F}$ is a final object and put $\mathcal{G} = i^{-1}\mathcal{F}$. It suffices to check that the canonical map $\mathcal{F} \to i_*i^{-1}\mathcal{F}$ is an isomorphism and it suffices to do this on the stalks. Let $x:X \to \mathcal{X}$ be a presentation and let $\xi$ be a geometric point of $X$. If $\xi$ lifts to a geometric point $\xi'$ of $U := X \times_{\mathcal{X}} \mathcal{U}$, then on the one hand $\mathcal{F}_{(X,x),\xi} = (j^{-1}\mathcal{F})_{(U,x|U),\xi'}$ is a one point set, while on the other hand the same is true for $(i_*i^{-1}\mathcal{F})_{(X,x),\xi} = (j^{-1}i_*i^{-1}\mathcal{F})_{(U,x|U),\xi'}$ since we already checked that $j^{-1}i_*$ maps all sheaves to the final object. If $\xi$ does not lift to $X \times_{\mathcal{X}} \mathcal{U}$, then it lifts to a geometric point $\xi'$ of $Z := X \times_{\mathcal{X}} \mathcal{Z}$, since $Z$ is the complement of $U$. So

$$(i_*i^{-1}\mathcal{F})_{(X,x),\xi} = (i^{-1}i_*i^{-1}\mathcal{F})_{(Z,x|Z),\xi'} \xrightarrow{\sim} (i^{-1}\mathcal{F})_{(Z,x|Z),\xi'} = \mathcal{F}_{(X,x),\xi}. \qquad \square$$

As a consequence of this proposition, if $i:\mathcal{Z} \to \mathcal{X}$ is a closed embedding, the induced morphism $\mathcal{Z}_{\text{ét}} \to \mathcal{X}_{\text{ét}}$ is a closed embedding.

Just as in SGA 4 VIII 6.5, we can now apply the theory of SGA 4 IV to obtain the usual functors $i^!$ for closed embeddings and $j_!$ for open embeddings, which satisfy the same properties as the corresponding constructions for schemes as given in [ibid.] or in [Mi, pp. 73–76].

## 10. Appendix: universal submersions

A property of morphisms between schemes that is of local nature for the smooth topology (see [L-MB, p. 33]) extends to a property of morphisms between algebraic stacks. If this property is also of topological nature, then one can check if it it is compatible with the formation of the topological space associated to a stack. We will treat the property 'being a universal submersion' in this subsection. This example is omitted in [L-MB].

Recall (EGA IV 15.7.8) that a morphism $f: X \to Y$ of schemes is called *submersive* if the associated map between topological spaces is surjective and a quotient map, i.e., $U \subset Y$ is open if $f^{-1}(U)$ is open. The map $f$ is *universally submersive* if it stays submersive after an arbitrary base change $Y' \to Y$.

**10.1. Proposition.** *Consider a commutative diagram in the category of schemes*

$$\begin{array}{ccccc} X'' & \xrightarrow{p'} & X' & \longrightarrow & X \\ & \searrow^{f''} & \downarrow & \square & \downarrow f \\ & & Y' & \xrightarrow[q]{} & Y \end{array}$$

*with $p'$ and $q$ surjective and smooth, and with the square cartesian. Then $f$ is universally submersive if and only if $f''$ is universally submersive.*

PROOF. The maps $p'$ and $q$ are smooth and therefore universally open. Being surjective, they are universally submersive. Furthermore, a composition of universal submersions is universally submersive. It follows that $f''$ is universally submersive if $f$ is. Conversely, suppose $f''$ is universally submersive. Then the composition of $f$ with the map $X'' \to X$ is universally submersive, and this implies that $f$ itself is universally submersive. $\square$

Hence we can extend the notion of being universally submersive to morphisms of algebraic spaces ([Kn, II]) and algebraic $S$-stacks ([L-MB, 4.14]). So by definition, a morphism $f: \mathcal{X} \to \mathcal{Y}$ of algebraic $S$-stacks is *universally submersive* if for some presentation $Y \to \mathcal{Y}$ and for some presentation $X' \to \mathcal{X} \times_{\mathcal{Y}} Y$, with $X'$ and $Y$ schemes, the induced morphism $X' \to Y$ is universally submersive. Note that this notion is stable under arbitrary base change in the 2-category of algebraic $S$-stacks.

Recall from [L-MB, §5] the construction of a functor $|\cdot|$ from the 2-category of algebraic $S$-stack to the category of topological spaces. For an algebraic $S$-stack $\mathcal{X}$, the open subsets of $|\mathcal{X}|$ correspond to the open substacks of $\mathcal{X}$.

**10.2. Proposition.** *Let $F: \mathcal{X} \to \mathcal{Y}$ be a morphism of algebraic $S$-stacks. For $F$ to be universally submersive, a necessary and sufficient condition is that for every morphism $G: \mathcal{Y}' \to \mathcal{Y}$ of algebraic $S$-stacks, when denoting by $F': \mathcal{X}' \to \mathcal{Y}'$ the morphism obtained by base-change the continuous map $|F'|: |\mathcal{X}'| \to |\mathcal{Y}'|$ is surjective and a quotient map.*

PROOF. Consider a 2-commutative diagram

(10.3)
$$\begin{array}{ccccc} X'' & \xrightarrow{h} & \mathscr{X}'' & \xrightarrow{g'} & \mathscr{X}' \\ & \searrow^{f''} & \downarrow & \square & \downarrow F' \\ & & Y' & \xrightarrow[g]{} & \mathscr{Y}' \end{array}$$

where $h$ and $g$ are presentations, where $X''$ and $Y'$ are schemes, and where the square is 2-cartesian.

Suppose $F$ is universally submersive; then so is $f''$. By [L-MB, 5.4ii] $|F'|$ is surjective. Let $U \subset |\mathscr{Y}'|$ and suppose $|F'|^{-1}U$ is open. Then $|f''|^{-1}|g|^{-1}U = |g'h|^{-1}|F'|^{-1}U$ is open and as $f''$ is submersive, $|g|^{-1}U$ is open. So by [L-MB, 5.6.1i] $U$ is open.

Now we prove the inverse implication. Take for $G\colon \mathscr{Y}' \to \mathscr{Y}$ a presentation with $\mathscr{Y}'$ a scheme. In (10.3) we will now suppose that $g$ and $g'$ are the identities (and $h$ is still a presentation of $\mathscr{X}'' = \mathscr{X}'$). Now $|F'|$ is by assumption surjective and a quotient map and by [ibid.] the same is true for $|h|$. So their composite $|f''|$ is also a surjective quotient map and hence $f''$ is submersive. If $\mathscr{Y}'' \to \mathscr{Y}$ is a morphism and we repeat this argument with $\mathscr{Y}'$ replaced by $Y' \times_{\mathscr{Y}} \mathscr{Y}''$, we see that $f''$ stays submersive after base change. But then $F$ is universally submersive by definition. □

## References for chapter 2